\documentclass[final]{siamltex}
\usepackage{amsmath}
\usepackage{amsfonts}
\usepackage{graphicx}

\newcommand*{\hD}{{\hat{D}}}
\newcommand*{\hP}{{\hat{P}}}
\newcommand*{\hU}{{\hat{U}}}
\newcommand*{\hu}{{\hat{u}}}
\newcommand*{\hv}{{\hat{v}}}

\newcommand*{\cC}{{\mathcal{C}}}
\newcommand*{\cL}{{\mathcal{L}}}

\newcommand*{\tL}{{\tilde{L}}}

\newcommand*{\RR}{{\mathbb{R}}}
\newcommand*{\CC}{{\mathbb{C}}}
\newcommand*{\ZZ}{{\mathbb{Z}}}
\newcommand*{\Vb}{{\mathbb{V}}}
\newcommand*{\Ub}{{\mathbb{U}}}

\newcommand*{\Gr}{{\mathrm{Gr}}}
\newcommand*{\GL}{{\mathrm{GL}}}
\newcommand*{\ib}{{\mathfrak i}}

\newcommand{\ra}{{\rightarrow}}
\newcommand{\I}{{I}}

\title{Computing stability of multi-dimensional 
travelling waves\thanks{VL is a postdoctoral fellow of the 
Fund of Scientific Research---Flanders (F.W.O.--Vlaanderen).
SJAM and JN were supported by EPSRC First Grant GR/S22134/01.
SJAM was also supported by Nuffield Foundation Newly Appointed 
Science Lecturers Grant SCI/180/97/129 and  
London Mathematical Society Scheme 2 Grant 2611.
JN was also supported by Australian Research Council grant DP0559083.}}

\author{Veerle Ledoux\thanks{Vakgroep 
Toegepaste Wiskunde en Informatica, Ghent University,
Krijgslaan 281-S9, B-9000 Gent, Belgium (Veerle.Ledoux@ugent.be)} 
\and Simon J.A. Malham\thanks{Maxwell Institute of Mathematical Sciences
and School of Mathematical and Computer Sciences 
Heriot-Watt University, Edinburgh EH14 4AS, UK (S.J.Malham@ma.hw.ac.uk)}
\and Jitse Niesen\thanks{School of Mathematics, 
University of Leeds, Leeds, LS2 9JT, UK (jitse@maths.leeds.ac.uk)}
\and Vera Th\"ummler\thanks{Fakult\"at f\"ur Mathematik, 
Universit\"at Bielefeld,
33501 Bielefeld, Germany (thuemmle@math.uni-bielefeld.de)}}

\begin{document}
\maketitle

\begin{abstract}
  We present a numerical method for computing the pure-point spectrum
  associated with the linear stability of multi-dimensional travelling
  fronts to parabolic nonlinear systems. Our method is based on the
  Evans function shooting approach.  Transverse to the direction of
  propagation we project the spectral equations onto a finite Fourier
  basis.  This generates a large, linear, one-dimensional system of
  equations for the longitudinal Fourier coefficients.  We construct
  the stable and unstable solution subspaces associated with the
  longitudinal far-field zero boundary conditions, retaining only the
  information required for matching, by integrating the Riccati
  equations associated with the underlying Grassmannian manifolds.
  The Evans function is then the matching condition measuring the
  linear dependence of the stable and unstable subspaces and thus
  determines eigenvalues. As a model application, we study the
  stability of two-dimensional wrinkled front solutions to a cubic
  autocatalysis model system. We compare our shooting approach with
  the continuous orthogonalization method of Humpherys and Zumbrun.
  We then also compare these with standard projection methods that
  directly project the spectral problem onto a finite
  multi-dimensional basis satisfying the boundary conditions.
\end{abstract}

\begin{keywords}
multi-dimensional stability, parabolic systems, Evans function
\end{keywords}

\begin{AMS}
65L15, 65L10
\end{AMS}

\markboth{Ledoux, Malham, Niesen and Th\"ummler}{Multi-dimensional stability}
\pagestyle{myheadings}

\thispagestyle{plain}

\section{Introduction}
Our goal in this paper is to present a
practical extension of the Evans function 
shooting method to computing the stability 
of \emph{multi-dimensional} travelling wave solutions 
to parabolic nonlinear systems. 
The problem is thus to solve the linear spectral equations 
for small perturbations of the travelling wave.  
Hence we need to determine the values of the spectral parameter 
for which solutions \emph{match} longitudinal far-field
zero boundary conditions and, in more than one spatial dimension,
the transverse boundary conditions. These are the \emph{eigenvalues}. 
There are two main approaches to solving such eigenvalue problems.
\begin{itemize}
\item \emph{Projection:} Project the spectral equations onto a finite dimensional
basis, which by construction satisfies the boundary conditions.
Then solve the resulting matrix eigenvalue problem.
\item \emph{Shooting:} Start with a specific value of the spectral parameter and 
the correct boundary conditions at one end. Shoot/integrate towards
the far end. Examine how close this solution is to 
satisfying the boundary conditions at the far end.
Re-adjust the spectral parameter accordingly.
Repeat the shooting procedure until the solution
\emph{matches} the far boundary conditions.
\end{itemize}
The projection method can be applied to travelling 
waves of any dimension. Its main disadvantage is that to
increase accuracy or include adaptation is costly and complicated
as we have to re-project onto the finer or adapted basis. Further,
from the set of eigenvalues that are produced by the resulting 
large algebraic eigenvalue problem, care must be taken to distinguish
two subsets as follows. In the limit of vanishing discretization scale,
one subset corresponds to isolated eigenvalues of finite multiplicity,
while the other subset converges to the absolute or essential
spectrum, depending on the boundary conditions (see Sandstede and
Scheel~\cite{SaSch}).
The latter subset can be detected by changing the truncated domain
size and observing which eigenvalues drift in the complex spectral
parameter plane. In contrast, for the shooting method,
it is much easier to fine-tune accuracy and adaptivity,  
and zeros of the matching condition are isolated eigenvalues of 
finite multiplicity (we do not have to worry about spurious
eigenvalues as described above). Its disadvantage though is that by 
its very description it is a one-dimensional method.

When implementing the shooting method, a Wronskian-like determinant
measuring the discrepancy in the boundary \emph{matching condition}
in the spectral problem is typically used to locate eigenvalues
in the complex spectral parameter plane. This is known as
the Evans function or miss-distance function---see 
Alexander, Gardner and Jones~\cite{AGJ}
or Greenberg and Marletta~\cite{GM}.
It is a function of the spectral parameter 
whose zeros correspond to eigenvalues, i.e.\ a \emph{match}. 
In practice, we match at a point roughly centred at the front.
Further for parabolic systems it is typically analytic
in the right-half complex parameter plane and the
argument principle allows us to globally determine 
the existence of zeros and therefore unstable eigenvalues
by integrating the Evans function along suitably chosen
contours.

The one-dimensional label of the Evans function has started to
be overturned and in particular Deng and Nii~\cite{DN},
Gesztesy, Latushkin and Makarov~\cite{GLM}, and
Gesztesy, Latushkin and Zumbrun~\cite{GLZ}
have extended the Evans function theoretically to multi-dimensions. 
Importantly from a computational
perspective, Humpherys and Zumbrun~\cite{HZ} 
outlined a computational method crucial to its 
numerical evaluation in the multi-dimensional scenario.
Their method overcame a restriction on the order
of the system that could be investigated numerically
using shooting, in particular using exterior product spaces. 
It is based on
continuous orhogonalization and is related to 
integrating along the underlying Grassmannian manifold
whilst evolving the coordinate representation for the
Grassmannian. 
More recently Ledoux, Malham and Th\"ummler~\cite{LMT} 
provided an alternative (though related) approach to that 
of Humpherys and Zumbrun. This approach is based on chosen 
coordinatizations of the Grassmannian manifold 
and integrating the resulting Riccati systems.

The main idea of this paper is to determine the
stability of multi-dimensional travelling waves by
essentially combining the two main approaches in the 
natural way: we project the spectral problem 
transversely and shoot longitudinally.
The basic details are as follows.
\begin{enumerate}
\item Transverse to the direction of propagation 
of the travelling wave we project onto 
a finite Fourier basis. This generates a large, 
linear, one-dimensional system of equations for 
the longitudinal Fourier coefficients.
\item We construct the stable and unstable solution 
subspaces associated with the longitudinal 
far-field zero boundary conditions by integrating
with given Grassmannian coordinatizations (chosen 
differently for each subspace). 
\item The Evans function is the determinant of the
matrix of vectors spanning both subspaces and measures
their linear dependence. 
We evaluate it at an intermediate longitudinal point.
\end{enumerate}
We compare the construction of the stable and unstable solution 
subspaces by integrating the spectral problem using the 
continuous orthogonalization method of Humpherys and Zumbrun.
We also compare these with standard projection methods 
that project the spectral problem onto 
a finite multi-dimensional basis satisfying the 
boundary conditions and  
solve the resulting large algebraic eigenvalue problem.

The parabolic nonlinear systems that we will consider here
are of the form
\begin{equation}\label{eq:sys}
\partial_tU=B\,\Delta U+c\,\partial_xU+F(U),
\end{equation}
on the cylindrical domain $\mathbb R\times\mathbb T$.
Here $U$ is a vector of component fields and $F(U)$
represents a nonlinear coupling reaction term. The
diagonal matrix $B$ encodes the positive constant diffusion 
coefficients.
We suppose we have Dirichlet boundary conditions in the
infinite longitudinal direction $x\in\mathbb R$ 
and periodic boundary conditions in the
transverse direction $y\in\mathbb T$. 
We assume also, that we are in a reference
frame travelling in the longitudinal direction 
with velocity $c$.
Though we assume only two spatial dimensions
our subsequent analysis in the rest of this paper
is straightforwardly extended to higher dimensional 
domains which are transversely periodic.
We also suppose we have a travelling wave solution 
$U=U_c(x,y)$ of velocity $c$ satisfying the boundary conditions, 
whose stability is the object of our scrutiny.

As a specific application we consider a cubic autocatalysis
reaction-diffusion system that admits wrinkled cellular
travelling fronts. They appear as solutions bifurcating 
from planar travelling wave solutions as the autocatalytic
diffusion parameter is increased. The initiation of this
first bifurcation is known as the cellular instability
and has been studied in some detail in the reaction kinetics
and combustion literature; see for example 
Horv\'ath, Petrov, Scott and Showalter~\cite{HPSS}
or Terman~\cite{T}. However our main interest here
is the stability of these cellular fronts themselves.
It has been shown by direct numerical
simulation of the reaction-diffusion system that the 
cellular fronts become unstable
as the autocatalyst diffusion parameter is 
increased further---see for example
Horv\'ath, Petrov, Scott and Showalter~\cite{HPSS}
and Malevanets, Careta and Kapral~\cite{MCK}.
We apply our multi-dimensional Evans function
shooting method to this problem and we reveal 
the explicit form of the spectrum associated 
with such fronts. We compare our shooting methods 
with standard direct projection methods obtaining 
high accuracy concurrency. 
We confirm the literature on these instabilities 
and establish the shooting method as an
effective and accurate alternative to 
standard direct projection methods.

Our paper is structured as follows.
In Section~\ref{sec:evans} we review 
Grassmannian and Stiefel manifolds.
We show how to construct the stable and
unstable solution subspaces of the spectral
problem using a given Grassmannian coordinatization
and define the Evans function using 
coordinatizations for each subspace.
We show how to project the spectral
problem onto a finite dimensional transverse
basis in detail in Section~\ref{sec:proj}.
A large portion of numerical effort is
devoted to accurately constructing the 
multi-dimensional fronts whose stability
we wish to study, and in Section~\ref{sec:comp}
we detail the technique we used to do this.
In Section~\ref{sec:auto} we outline our
main application---the cubic autocatalysis
problem---and review the known travelling
wave solutions and stability properties. 
We then bring to bear our full range of numerical
techniques to the cubic autocatalysis
problem in Section~\ref{sec:numerics},
comparing their accuracy.
In Section~\ref{sec:conclu} we provide 
some concluding remarks and outline future
directions of investigation.

\section{Review: Grassmann flows and spectral matching}\label{sec:evans}

\subsection{Grassmannian coordinate representations}
A $m$-frame is a $m$-tuple of $m\leq n$ linearly independent
vectors in $\CC^n$. The \emph{Stiefel manifold} $\Vb(n,m)$ of $m$-frames 
is the open subset of $\CC^{n\times m}$ of all $m$-frames 
centred at the origin. Hence any element $Y\in\Vb(n,m)$
can be represented by an $n\times m$ matrix of rank $m$:
\begin{equation}\label{origrep}
Y=\begin{pmatrix} Y_{11} & \cdots & Y_{1m} \\
\vdots & & \vdots \\
Y_{n1} & \cdots & Y_{nm}
\end{pmatrix}.
\end{equation}
The set of $m$ dimensional subspaces of $\CC^n$
forms a complex manifold $\mathrm{Gr}(n,m)$ called 
the \emph{Grassmann manifold} of $m$-planes in $\CC^n$;
it is compact and connected (see Steenrod~\cite[p.~35]{St} 
or Griffiths and Harris~\cite[p.~193]{GH}).
There is a quotient map from $\Vb(n,m)$ to $\Gr(n,m)$,
sending each $m$-frame centred at the origin to the $m$-plane
it spans---see Milnor-Stasheff~\cite[p.~56]{MS}.
Any $m$-plane in $\CC^n$ can be represented by
an $n\times m$ matrix of rank $m$, like $Y$ above.
However any two such matrices $Y$ and $Y'$ related 
by a rank $m$ transformation $u\in\GL(k)$, so that $Y'=Yu$,
will represent the same $m$-plane. Hence we can cover 
the Grassmann manifold $\Gr(n,m)$ by coordinate patches
$\Ub_\ib$, labelled with multi-index 
$\ib=\{i_1,\ldots,i_m\}\subset\{1,\ldots,n\}$,
where $\Ub_\ib$ is the set of $m$-planes $Y\in\Gr(n,m)$,
which have a matrix representation $y_{\ib^\circ}\in\CC^{n\times m}$
whose $m\times m$ submatrix, designated by the rows $\ib$, is 
the identity matrix. Each such patch 
is an open dense subset of $\Gr(n,m)$ isomorphic to $\CC^{(n-m)m}$.
For example, if $\ib=\{1,\ldots,m\}$, then $\Ub_{\{1,\ldots,m\}}$
can be uniquely represented by a matrix of the form
\begin{equation}\label{eq:patch1}
y_{\ib^\circ}=\begin{pmatrix} \I_m \\
                            \hat y 
            \end{pmatrix}
\end{equation}
with $\ib^\circ=\{m+1,\ldots,n\}$ and $\hat y\in\CC^{(n-m)m}$. 
For each $\ib$, there is a bijective map 
$\varphi_{\ib}\colon\Ub_\ib\ra\CC^{(n-m)m}$ given by 
$\varphi_{\ib}\colon y_{\ib^\circ}\mapsto \hat y$,
representing the local coordinate chart for the coordinate patch $\Ub_\ib$.
For more details see Griffiths and Harris~\cite[p.193--4]{GH}.

\subsection{Riccati flow}
Consider a non-autonomous linear vector field defined
on the Stiefel manifold of the form $V(x,Y)=A(x)\,Y$,
where $Y\in\Vb(n,m)$ and $A\colon\mathbb R\ra\mathfrak{gl}(n)$,
the general linear algebra of rank $n$ matrices.
Following the exposition in Ledoux, Malham and Th\"ummler~\cite{LMT},
fix a coordinate patch $\Ub_\ib$ for $\Gr(n,m)$ for some $\ib$,
and decompose $Y\in\Vb(n,m)$ into
\begin{equation*}
Y=y_{\ib^\circ} u.
\end{equation*}
If we substitute this form into the ordinary differential
system $Y'=V(x,Y)$ we obtain:
\begin{equation*}
y_{\ib^\circ}' u+y_{\ib^\circ} u'=(A_\ib+A_{\ib^\circ}y_{\ib^\circ})\,u,
\end{equation*}
where $A_\ib$ denotes the submatrix obtained by restricting the matrix
$A(x)$ to its $\ib$th columns. If we project this last system onto
its $\ib^\circ$th and $\ib$th rows, respectively, we generate 
the following pair of equations for the coordinate chart variables
$\hat y=\varphi_\ib\circ y_{\ib^\circ}$ and transformations $u$,
respectively:
\begin{subequations}
\begin{align}\label{eq:riccati}
\hat y'&=c+d\hat y-\hat ya-\hat yb\hat y\\ 
\intertext{and}
u'&=(a+b\hat y)\,u.
\end{align}
\end{subequations}
Here $a$, $b$, $c$ and $d$ denote
the $\ib\times\ib$, $\ib\times \ib^\circ$, $\ib^\circ\times \ib$
and $\ib^\circ\times \ib^\circ$ submatrices of $A$, respectively.

We observe that the flow on the Stiefel manifold $\Vb(n,m)$ generated
by the linear vector field $V(x,Y)$ can be decomposed into an 
\emph{independent Riccati flow} in the coordinate chart variables of a fixed
coordinate patch $\Ub_\ib$ of $\Gr(n,m)$, and a 
\emph{slaved linear flow} of transformations in $\GL(m)$. 
Hence we can reconstruct the flow on the Stiefel manifold
generated by the linear vector field $V(x,Y)$, by simply
solving the Riccati equation above, and then if required,
subsequently integrating the linear equation for the 
transformations. In general, solutions to the
Riccati equation~\eqref{eq:riccati} can blow up.
However, this is simply the manifestation of
a poorly chosen local representative coordinate patch.
This can be gleaned from the fact that 
simultaneously the determinant of $u$ becomes zero---afterall
the linear Stiefel flow, with globally smooth coefficients, 
does not blow up. To resolve this, we simply change patch
and keep a record of the determinant of the patch swapping 
transformation. Further details and 
strategies can be found in Ledoux, Malham and Th\"ummler~\cite{LMT}.

\subsection{Drury--Oja flow}
The continuous orthogonalization method of 
Humpherys and Zumbrun~\cite{HZ} can be thought
of as the flow on the Grassmann manifold 
corresponding to the linear vector field $V$,
with the coordinatization evolving according to a given
unitary flow (see Ledoux, Malham and Th\"ummler~\cite{LMT}). 
Indeed their method is specified by a Drury--Oja flow
for the $n\times m$ orthogonal matrix $Q$, together 
with the flow for the determinant of the matrix $R$,
in the $QR$-decomposition of $Y$:
\begin{subequations}\label{eq:DO}
\begin{align}
Q'&=(\I_n-QQ^\dag)A(x)Q,\\
(\det R)'&=\mathrm{Tr}\bigl(Q^\dag A(x)Q\bigr)(\det R).
\end{align}
\end{subequations}
Here $Q^\dag$ denotes the Hermitian transpose of $Q$.
In practice the determinant of $R$ is exponentially rescaled---see 
Humpherys and Zumbrun~\cite{HZ}. We also did not find it necessary to apply
any of the stabilization techniques suggested therein.

\subsection{Spectral problem}
Consider the linear spectral problem on $\RR$ with 
spectral parameter $\lambda$ in 
standard first order form with coefficient
matrix $A(x;\lambda)\in\CC^{n\times n}$:
\begin{equation}\label{eq:linearevans} 
Y'=A(x;\lambda)\,Y.
\end{equation}
We assume there exists a subdomain $\Omega\subset\CC$
containing the right-half complex plane, that does not
intersect the essential spectrum. 
For $\lambda\in\Omega$, we know that there
exists exponential dichotomies on $\mathbb R^-$ and 
$\mathbb R^+$ with the same Morse index $m$ in each case
(see Sandstede~\cite{Sand}). Hence for $\lambda\in\Omega$,
let $Y^-(x;\lambda)\in\Vb(n,m)$ denote 
the matrix whose columns are solutions to~\eqref{eq:linearevans} 
which span the unstable subspace of solutions 
decaying exponentially to zero as $x\ra-\infty$, and 
let $Y^+(x;\lambda)\in\Vb(n,n-m)$ denote the matrix 
whose columns are the solutions which span the stable
subspace of solutions decaying exponentially 
to zero as $x\ra+\infty$.
The values of spectral parameter $\lambda$ for which the 
columns of $Y^-$ and columns of $Y^+$ are linearly dependent
on $\RR$ are pure-point eigenvalues.

\subsection{Evans determinant}\label{sec:evans-determinant}
The Evans function $D(\lambda)$ is the measure of the
linear dependence between the two basis sets $Y^-$ and $Y^+$:
\begin{equation}\label{defevans}
D(\lambda)\equiv 
\mathrm{e}^{-\int_0^x\mathrm{Tr}A(\xi;\lambda)\,\mathrm{d}\xi}\,
\mathrm{det}\bigl(Y^-(x;\lambda)\,\, Y^+(x;\lambda)\bigr),
\end{equation}
(see Alexander, Gardner and Jones~\cite{AGJ} or Sandstede~\cite{Sand}
for more details). Note that for the decompositions 
$Y^-=y_{\ib_-^\circ}u_{-}$ and $Y^+=y_{\ib_+^\circ}u_{+}$,
assuming $\det\,u_{-}\neq0$ and $\det\,u_{+}\neq0$, we have
\begin{equation*}
\det\bigl(Y^-\, Y^+\bigr)
=\det\begin{pmatrix} y_{\ib_-^\circ} & y_{\ib_+^\circ} \end{pmatrix}
\cdot\det\,u_{-}\cdot\det\,u_{+}.
\end{equation*}
Let $Y_0^-(\lambda)$ denote the $n\times m$ matrix whose 
columns are the $m$ eigenvectors of $A(-\infty;\lambda)$
corresponding to eigenvalues with a positive real part.
Analogously let $Y_0^+(\lambda)$ denote the $n\times(n-m)$ 
matrix whose columns are the $(n-m)$ eigenvectors of $A(+\infty;\lambda)$
corresponding to eigenvalues with a negative real part.

Suppose we fix a coordinate patch for $\Gr(n,m)$ 
labelled by $\ib_-$ (which has cardinality $m$). 
We integrate the corresponding
Riccati differential problem~\eqref{eq:riccati} for
$\hat y_-=\varphi_{\ib_-}\circ y_{\ib_-^\circ}$ from $x=-\infty$
to $x=x_\ast$. The initial data
$\hat y_-(-\infty;\lambda)$ is generated  
by postmultiplying the $\ib_-^\circ\times\{1,\ldots,m\}$
submatrix of $Y_0^-(\lambda)$ by the inverse of the $\ib_-\times\{1,\ldots,m\}$
submatrix of $Y_0^-(\lambda)$ (i.e.\ perform the simple decomposition into 
$y_{\ib_-^\circ}$ and $u_{-}$ for $Y_0^-(\lambda)$, and 
use the $\ib_-^\circ\times\{1,\ldots,m\}$ block of 
$y_{\ib_-^\circ}$ as the initial data for $\hat y_-$).
Similarly we fix a patch for $\Gr(n,n-m)$ labelled by $\ib_+$
(with cardinality $(n-m)$) and integrate the 
Riccati differential problem~\eqref{eq:riccati} for
$\hat y_+=\varphi_{\ib_+}\circ y_{\ib_+^\circ}$ from $x=+\infty$
to $x=x_\ast$. We perform the analogous decomposition on
$Y^+_0(\lambda)$ with index $\ib_+$ to generate the initial data 
$\hat y_+(+\infty;\lambda)$. If the solutions to both 
Riccati problems do not blow up, so that 
$\det u_{-}(x;\lambda)\neq0$ for all $x\in(-\infty,x_\ast]$ 
and $\det u_{+}(x;\lambda)\neq0$ for all $x\in[x_\ast,\infty)$, 
then we define the modified Evans function, 
which is analytic in $\lambda$, by
\begin{equation}\label{newdefevans}
D(\lambda;x_\ast)\equiv 
\det\begin{pmatrix} 
              y_{\ib_-^\circ}(x_\ast;\lambda) &  y_{\ib_+^\circ}(x_\ast;\lambda) 
\end{pmatrix},
\end{equation}
where $x_\ast\in\RR$ is the matching point. 
In our practical implementation in Section~\ref{sec:auto},
we took the patches labelled by $\ib_-=\{1,\ldots,m\}$ 
and $\ib_+=\{m+1,\ldots,n\}$ for the left and right problems,
respectively. We used the generally accepted rule of thumb
for Evans function calculations: 
to match at a point roughly centred at the front.
In our example application we took $x_*=50$, and did \emph{not} observe 
singularities in $\hat y^-$ or $\hat y^+$.
However generally, varying $x_\ast$ may introduce singularities for
$\hat y^-$ or $\hat y^+$ in their respective domains 
$(-\infty,x_\ast]$ and $[x_\ast,\infty)$. Indeed this was observed 
in several examples considered in Ledoux, Malham and Th\"ummler~\cite{LMT}.
There, a robust remedy is provided, which on numerical evidence,
allows matching anywhere in the computational domain (and which
in future work we intend to apply to our context here). 

If we use the Humpherys and Zumbrun continuous orthogonalization
approach then we define the Evans function to be
\begin{equation}\label{eq:HZEvans}
D_{\text{HZ}}(\lambda;x_\ast)\equiv\det \bigl(Q^-(x_\ast;\lambda)~Q^+(x_\ast;\lambda)\bigr)
\cdot\det R^-(x_\ast;\lambda)\cdot\det R^+(x_\ast;\lambda).
\end{equation}
This is analytic in $\lambda$ when we include the $\det R^\pm$ terms.

\subsection{Angle between subspaces}
We can also measure the angle between the unstable 
and stable linear subspaces $\Vb(n,m)$ and $\Vb(n,n-m)$,
see Bj\"orck and Golub~\cite{BG}. When the columns of the matrices
$Y^-(x;\lambda)$ and $Y^+(x;\lambda)$ which span 
$\Vb(n,m)$ and $\Vb(n,n-m)$, respectively,  
are nearly linearly dependent the angle will be small.
It is defined as follows---without loss of generality 
we assume $m\geq n/2$. Let $Q^-$ and $Q^+$ be unitary bases for 
$\Vb(n,m)$ and $\Vb(n,n-m)$. They are specified by the 
$QR$-decompositions $Y^-=Q^-R^-$ and $Y^+=Q^+R^+$.
The cosine of the smallest angle $\theta$ 
between $\Vb(n,m)$ and $\Vb(n,n-m)$ is given by
the largest singular value of the matrix $(Q^-)^\dag Q^+$.
However when $\theta$ is small it is more accurate 
to compute $\sin\theta$, which 
is given by the smallest singular value of 
\begin{equation*}
\big(\I_n-Q^-(Q^-)^\dag\bigr)Q^+.
\end{equation*}
Hence if the singular values of this matrix
are $\sigma_i$, $i=1,\ldots,n-m$, define
\begin{equation*}
\theta(\lambda;x_\ast)\equiv\arcsin\min_{i=1,\ldots,n-k}\sigma_i,
\end{equation*}
where again, $x_\ast$ is the matching point.
Eigenvalues correspond to zeros of $\theta(\lambda;x_\ast)$.

\section{Projection onto a finite transverse Fourier basis}\label{sec:proj}
Consider the parabolic nonlinear system~\eqref{eq:sys}
posed on the cylindrical domain $\mathbb R\times\mathbb T$.
Recall that we assume Dirichlet boundary conditions in the
infinite longitudinal direction $x\in\mathbb R$ 
and periodic boundary conditions in the
transverse direction $y\in\mathbb T$ (here $\mathbb T$ 
is the one-torus with period/circumference $L$). 
Suppose we have a travelling wave solution $U=U_c(x,y)$
of velocity $c$ satisfying the boundary conditions.
The stability of this solution is determined 
by the spectrum of the linear differential operator
\begin{equation}
\label{eq:call}
\mathcal L=B\Delta+c\,\partial_x+\mathrm{D}F(U_c),
\end{equation}
where $\mathrm{D}F$ is the Jacobian of the interaction term.
Hence the eigenvalue problem is
\begin{equation}\label{eq:evalprob}
B\Delta U+c\,\partial_x U+\mathrm{D}F(U_c)U=\lambda U.
\end{equation}
We wish to define and compute the Evans function for this spectral problem.
The idea is to perform a Fourier decomposition in the 
transverse direction. This transforms the problem
to a one-dimensional problem, for which we can apply the Evans function
shooting approach. Hence suppose that we can write
\begin{equation}
\label{hatu}
U(x,y)=\sum_{k=-\infty}^\infty\hU_k(x)\,\mathrm{e}^{iky/\tL},
\end{equation}
where $\hU_k\in\mathbb C^N$ are the Fourier coefficients
(with $N$ denoting the number of components of~$U$)
and $\tL = 2\pi L$---this implicitly restricts our attention to
perturbations~$U$ with period~$L$. 
Similarly, we expand $\mathrm{D}F(U_c)$ as 
\begin{equation}
\label{hatd}
\mathrm{D}F(U_c(x,y))=\sum_{k=-\infty}^\infty\hD_k(x)\,\mathrm{e}^{iky/\tL},
\end{equation}
where the Fourier coefficients $\hD_k\in\mathbb C^{N\times N}$.
Substituting our Fourier expansions for $U$ and $\mathrm{D}F(U_c)$ 
in~\eqref{hatu} and~\eqref{hatd} in the 
eigenvalue problem~\eqref{eq:evalprob} yields
\begin{multline*}
B\,(\partial_{xx}+\partial_{yy})\sum_{k=-\infty}^\infty \hU_k
\,\mathrm{e}^{iky/\tL}+c\,\partial_x \sum_{k=-\infty}^\infty \hU_k\,\mathrm{e}^{iky/\tL}\\ 
+\sum_{\ell=-\infty}^\infty \hD_\ell \,\mathrm{e}^{i\ell y/\tL}
\sum_{k=-\infty}^\infty \hU_k \,\mathrm{e}^{iky/\tL}
=\lambda \sum_{k=-\infty}^\infty \hU_k \,\mathrm{e}^{iky/\tL}.
\end{multline*}
Reordering the double sum 
\begin{equation*}
\sum_{\ell=-\infty}^\infty \hD_\ell \,\mathrm{e}^{i\ell y/\tL}
 \sum_{k=-\infty}^\infty \hU_k \,\mathrm{e}^{iky/\tL} 
=\sum_{k=-\infty}^\infty \sum_{\nu=-\infty}^\infty \hD_{k-\nu} 
\hU_\nu \,\mathrm{e}^{iky/\tL},
\end{equation*}
yields that for all $k\in\ZZ$, we must have
\begin{equation*}
\partial_{xx}\hU_k-(k/\tL)^2\hU_k 
+ c\,B^{-1}\partial_x\hU_k 
+\sum_{\nu=-\infty}^\infty B^{-1}\hD_{k-\nu}\hU_\nu=\lambda B^{-1}\hU_k.
\end{equation*}
In practical implementation, we consider only 
the first~$K$ Fourier modes. Hence we replace the above
infinite system of equations by following finite system
of equations which we now also express in first order form: 
\begin{subequations}\label{eq:bigsys}
\begin{align}
\partial_x\hU_k&=\hP_k, \\
\partial_x\hP_k&=\lambda B^{-1}\hU_k 
+(k/\tL)^2\hU_k-c\, B^{-1}\hP_k
-\sum_{\nu=-K}^K B^{-1}\hD_{k-\nu}\hU_\nu,
\end{align}
\end{subequations}
for $k=-K,-K+1,\ldots,K$.
This is a large system of ordinary differential equations
involving the spectral parameter~$\lambda$, so
we can apply the standard Evans function approach to it. 
In matrix form this is
\begin{equation}\label{eq:bigsys2}
\partial_x \begin{pmatrix} \hU \\ \hP \end{pmatrix}
= \begin{pmatrix} O_{N(2K+1)} & I_{N(2K+1)} \\ A_3(x;\lambda) & A_4 \end{pmatrix}
\begin{pmatrix} \hU \\ \hP \end{pmatrix},
\end{equation}
where $\hU$ is the vector of $\mathbb C^N$-valued components
$\hU_{-K}$, $\hU_{-K+1}$, \ldots, $\hU_K$ and $\hP$
is the vector of $\mathbb C^N$-valued components
$\hP_{-K}$, $\hP_{-K+1}$, \ldots, $\hP_K$.
The lower right block matrix is $A_4=-c\,B^{-1}\otimes I_{2K+1}$.
The lower left block matrix is 
given by $A_3=\tilde A_3(\lambda)+\hat A_3(x)$ 
where if we set $E_k(\lambda)\equiv\lambda B^{-1}+(k/\tL)^2 I_N$ then
\begin{equation*}
\tilde A_3(\lambda)=
\begin{pmatrix}
E_{-K}(\lambda) & O & \cdots & O \\
O & E_{-K+1}(\lambda) & \cdots & O \\
\vdots & \vdots & \ddots & \vdots \\
O      & \cdots & O & E_{+K}(\lambda) 
\end{pmatrix}
\end{equation*}
and
\begin{equation*}
\hat A_3(x) = -B^{-1} \otimes \begin{pmatrix}
\hD_0    & \hD_{-1} & \cdots & \hD_{-2K} \\
\hD_{1}  & \hD_0 & \cdots & \hD_{-2K+1} \\
\vdots   & \vdots & \ddots & \vdots \\
\hD_{2K}  &\hD_{2K-1} & \cdots & \hD_0 
\end{pmatrix}.
\end{equation*} 
A similar Galerkin approximation for nonplanar fronts can be found in 
Gesztesy, Latushkin and Zumbrun~\cite[p.~28,37]{GLZ}.

\section{The freezing method}
\label{sec:comp}
To study the long-time behaviour of the wrinkled front
solution we use the method of freezing travelling waves described in
Beyn and Th\"ummler~\cite{BT} and analyzed in Th\"ummler~\cite{Th}. 
This method transports the problem into a moving frame whose velocity 
is computed during the computation. For the travelling waves whose
stability is our concern here, the freezing method has
two advantages over direct simulation: Firstly, we can do long-time
computations of time dependent (for example oscillating) solutions.
This is because the method picks out the correct speed for the 
moving frame in which the wave is roughly stationary (in particular
it does not drift out of the chosen domain).
Secondly, the freezing method can be used to obtain a suitably accurate 
guess for the waveform, which can then be used to initiate a Newton solver
that computes the travelling wave as a stationary solution in the
co-moving frame.   
For completeness we will briefly review the idea
of the freezing method here. Consider a parabolic nonlinear system of
form~\eqref{eq:sys} in a stationary frame of reference:
\begin{equation}\label{eq:pde}
\partial_tU=B\,\Delta U+F(U).
\end{equation}
If we substitute the ansatz $U(x,y,t)=V(x-\gamma(t),y,t)$ 
into this equation, where $\gamma(t)$ is the unknown reference position 
for a wave travelling in the longitudinal $x$-direction, we
generate the partial differential equation 
for $V$ and $\gamma'$:
\begin{subequations}\label{eq:pdae}
\begin{gather}
\partial_t V=B\,\Delta V+\gamma'(t)\partial_x V+F(V).\\ 
\intertext{In order to compensate for the additional degree of freedom
which has been introduced by the unknown position $\gamma$, we specify the
following additional algebraic constraint}
0=\int_{\RR\times\mathbb T} \bigl(\partial_x \hat V(x,y,t)\bigr)^{\mathrm{T}} 
     \bigl(\hat V(x,y,t)-V(x,y,t)\bigr)\,\mathrm{d}x\,\mathrm{d}y. 
\end{gather}
\end{subequations}
This constraint selects a unique solution from the one-parameter
group of possible solutions generated by longitudinal translation.
It is a phase condition that normalizes the system with respect
to a given template function---for example a suitable choice 
is to set $\hat V(x,y,t)$ to be the initial data $V(x,y,0)$.
Hence we have generated a partial differential algebraic 
system~\eqref{eq:pdae}, for the unknowns $V(x,y,t)$ and $\zeta(t)=\gamma'(t)$,
with initial data $V(x,y,0)$ and $\zeta(0)$ within the
attractive set of the stable travelling wave and its velocity.

A stationary solution $(V,\zeta)=(V_c,c)$ of the 
partial differential algebraic system~\eqref{eq:pdae}
represents a travelling wave solution $U_c(x,y,t)=V_c(x-ct,y)$ 
with fixed speed $c$ to the parabolic nonlinear problem~\eqref{eq:pde}.
If the travelling wave $U_c$ is stable then the stationary solution
$(V_c,c)$ is also stable and the time evolution of
partial differential algebraic system~\eqref{eq:pdae} 
will converge to $(V_c,c)$.
Thus the numerical method consists of
solving the partial differential algebraic system~\eqref{eq:pdae} 
until the solution $(V,\zeta)$ has reached the
equilibrium $(V_c,c)$. Alternatively one can also solve
the partial differential algebraic system~\eqref{eq:pdae} 
until some time $t=T$ when the solution is deemed sufficiently 
close to equilibrium and then use a direct method (e.g.\ Newton's
method) to solve the stationary equation 
\begin{subequations}\label{eq:equilib}
\begin{align}
 0 &=B\,\Delta v+\theta \partial_x v+F(v), \\
 0 &=\int_{\RR \times \mathbb T} \bigl(\partial_x \hat V(x,y)\bigr)^{\mathrm{T}} 
\bigl(\hat V(x,y)-v(x,y)\bigr)\,\mathrm{d}x\,\mathrm{d}y, 
\end{align}
\end{subequations}
for $(v,\theta)$ starting with initial values $v_0(x,y)=V(x,y,T)$ 
and $\theta_0=\zeta(T)$.

\section{Review: autocatalytic system}\label{sec:auto}
As our model application we consider the cubic autocatalysis system
\begin{subequations}
\label{cubicautoeqns}
\begin{align}
u_t&=\Delta u+c\partial_xu-uv^2,\\ 
v_t&=\delta\Delta v+c\partial_xv+uv^2.
\end{align}
\end{subequations}
Here $u(x,y,t)$ is the concentration of the 
reactant and $v(x,y,t)$ is the concentration 
of the autocatalyst, in the infinitely extended
medium $(x,y)\in\mathbb R\times\mathbb T$.
In the far field, we suppose $(u,v)$
approaches the stable homogenenous steady state $(0,1)$
as $x\rightarrow-\infty$, and the unstable homogeneous 
steady state $(1,0)$ as $x\rightarrow+\infty$. The
parameter $\delta$ is the ratio of the diffusivity
of the autocatalyst to that of the reactant. 
This system is is globally well-posed for smooth initial data,
and any finite $\delta>0$.

From Billingham and Needham~\cite{BN} we know that, in one spatial
dimension, a unique heteroclinic connection between 
the unstable and stable homogeneous steady states 
exists for wavespeeds $c\geq c_{\min}$. The
unique travelling wave for $c=c_{\min}$ converges
exponentially to the homogeneous steady states.
They are readily constructed by shooting 
(see for example Balmforth, Craster and Malham~\cite{BCM}) 
and represent the planar front solution referred to hereafter. 

It is well known that these planar travelling wave solutions
become unstable to transverse perturbations when $\delta$ is 
greater than a critical ratio $\delta_{\mathrm{cr}}\approx2.3$
(from Malevanets, Careta and Kapral~\cite[p.~4733]{MCK}),
and our transverse domain is sufficiently large.
This is known as the cellular instability---see for example
Terman~\cite{T}, Horv\'ath, Petrov, Scott and Showalter~\cite{HPSS}
and Malevanets, Careta and Kapral~\cite{MCK}.
Calculating the spectrum that reveals the mechanism
of the instability is also straightforward.
If we consider small perturbations
about the underlying planar travelling wave solution of the form
$\hat U(x)\mathrm{e}^{\lambda t+\mathrm{i}(k/\tilde L)y}$,
for any $k\in\mathbb Z$,
then the associated linear operator $\cL$ whose spectrum
we wish to determine has the form 
\begin{equation}
\label{eq:cL1d}
\cL=B\bigl(\partial_{xx}-(k/\tilde L)^2 I\bigr)+c\,\partial_{x}+\mathrm{D}F(U_c).
\end{equation}
The spectrum of $\cL$ consists of the pure-point spectrum,
i.e.\ isolated eigenvalues of finite multiplicity,
together with the essential spectrum. By a result of
Henry~\cite{H}, if we consider $\cL$ as an operator on
the space of square integrable functions,
we know that the essential spectrum is 
contained within the parabolic curves of the continuous
spectrum in left-half complex plane. 
The continuous spectrum for this model problem is 
given for all $\mu\in\mathbb R$ and each fixed $k\in\mathbb Z$, 
by the locus in the complex $\lambda$-plane of the curves
\begin{align*}
\lambda &=-\delta (k/\tilde L)^2+\mathrm{i}c\mu-\delta\mu^2,\\
\lambda &=-1-\delta (k/\tilde L)^2+\mathrm{i}c\mu-\delta\mu^2,\\
\lambda &=-(k/\tilde L)^2+\mathrm{i}c\mu-\mu^2,
\end{align*}
where the last curve appears with multiplicity two.
Hence the continuous spectrum touches the imaginary axis 
at the origin when $k$ is zero. 
The origin is a simple eigenvalue for $\cL$ due to 
translational invariance of the underlying
travelling wave $U_c(x)$ in the longitudinal direction.
Thus, even though the essential spectrum does not 
affect linear stability, we cannot automatically deduce 
asymptotic stability unless we consider 
a gently weighted space of square integrable functions, 
which will shift the essential spectrum to the left. 

\begin{figure}
  \begin{center}
    \includegraphics[width=\linewidth]{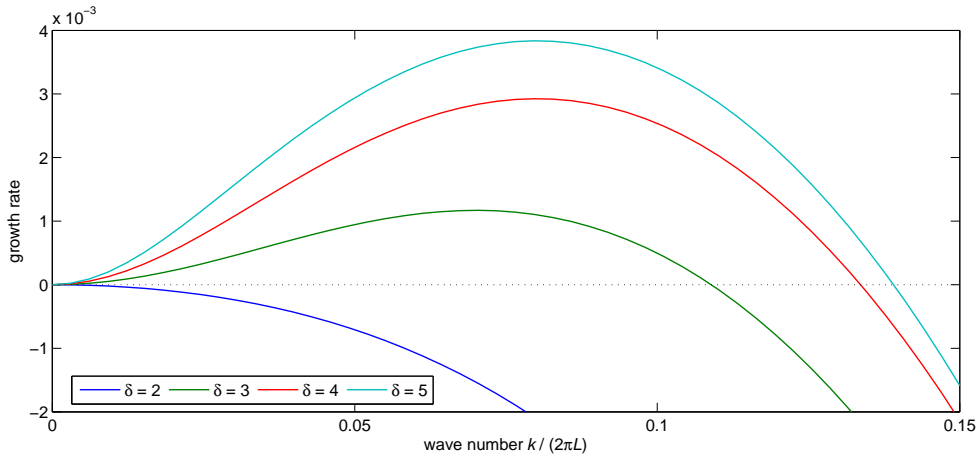} 
  \end{center}
  \caption{Dispersion relations for the planar front for different
    values of~$\delta$. For $\delta = 2$, all perturbations decay and
    the planar front is linearly stable, but for the other values
    of~$\delta$ there is a growing interval $(0,k_\delta)$ of unstable 
    transverse wavenumbers.} 
  \label{dispersion}
\end{figure}

We now notice that the spectral problem for the linear
differential operator $\cL$ is fourth order and one-dimensional, 
the transverse perturbation information is encoded though the 
transverse wavenumber parameter $k$. Hence we can employ the 
Evans function shooting approach, in particular in 
Figure~\ref{dispersion} we plot the dispersion relation
for planar fronts corresponding to different values of
the physical parameter $\delta$. The dispersion relation
gives the relationship between the real part of
the exponential growth factor of small perturbations,
$\mathrm{Re}(\lambda)$, and transverse wave number 
$k$ (in fact to make comparisons later easier we
have used $k/\tilde L$ on the abscissa). 
To construct the graph shown, we locate zeros
of the Evans function on the real $\lambda$ axis and 
follow them as the transverse wave number $k$ is varied.
We see in Figure~\ref{dispersion} that the
interval $(0,k_\delta)$ of unstable 
transverse wavenumbers $k$ increases as $\delta$ 
increases (at least up to $\delta=5$). The wavenumber where the
maximum of $\mathrm{Re}(\lambda)$ occurs has largest growth 
and is thus typically the transverse mode expressed 
in the dynamics. 

However our real concern here is the stability of the 
wrinkled cellular fronts themselves as we increase the
diffusion ratio parameter $\delta$ beyond $\delta_{\mathrm{cr}}$.
We know the planar waves are unstable in this regime
and the question is whether the bifurcating wrinkled
cellular fronts are stable for values of $\delta$ just
beyond $\delta_{\mathrm{cr}}$ but do become unstable at
another critical value of $\delta$. Using direct
simulations of the reaction-diffusion system
Malevanets, Careta and Kapral~\cite{MCK} demonstrated that
when $\delta=5$ the transverse structure of the propagating
fronts was very complex with spatio-temporal dynamics of 
its own. Our goal is to elicit in practice the initial
mechanisms that precipitate the behaviour they observe
through a careful study of the spectra of the underlying
multi-dimensional travelling fronts.

\section{Numerics}\label{sec:numerics}

\subsection{The method}\label{sec:method}
We will compute the spectrum of the two-dimensional travelling front
by computing the Evans function associated with the linear
operator~$\mathcal{L}$ defined in~\eqref{eq:call}. First, we compute
the front with the freezing method as explained in
Section~\ref{sec:comp}. Then as in Section~\ref{sec:proj}, we
project the problem onto a finite number
of transverse Fourier modes $k=-K,\ldots,K$. This generates the
large linear system~\eqref{eq:bigsys2} which is of 
the standard linear form~\eqref{eq:linearevans} where the coefficient
matrix is given by 
\begin{equation*}
A(x;\lambda)=\begin{pmatrix} O_{N(2K+1)} & I_{N(2K+1)} \\ 
                             A_3(x;\lambda) & A_4 \end{pmatrix},
\end{equation*}
where the matrices $A_3$ and $A_4$ are explicitly 
given at the end of Section~\ref{sec:proj}. This coefficient matrix 
depends on both the travelling wave solution~$U_c$ and the spectral
parameter~$\lambda$. It also depends on the nonlinear reaction term.
For our autocatalytic problem the Jacobian $\mathrm{D}F$ of the
nonlinear reaction term and its Fourier transform $\hD$ are given by
\begin{equation*}
\mathrm{D}F(U)=
\begin{pmatrix} 
-v^2 & -2uv \\ 
v^2 & 2uv 
\end{pmatrix}
\qquad\text{and}\qquad 
\hD = 
\begin{pmatrix} 
-\hv*\hv & -2\hu*\hv \\
\hv*\hv & 2\hu*\hv
\end{pmatrix},
\end{equation*}
where $\hu$ and $\hv$ are the Fourier transforms of~$u$ and~$v$
respectively, and $*$ denotes the convolution defined by 
$(\hu*\hv)_k=\sum_{\ell=-\infty}^\infty \hu_\ell \hv_{k-\ell}$. 

As described in Section~\ref{sec:evans-determinant},
we construct the unstable subspace $Y_0^-(\lambda)$ of~\eqref{eq:bigsys2}
at $x=-\infty$, as the matrix whose columns are the eigenvectors of 
$A(-\infty;\lambda)$ corresponding to eigenvalues with a
positive real part. Similarly we construct the stable subspace
$Y_0^+(\lambda)$ at $x=+\infty$ from the eigenvectors of 
$A(+\infty;\lambda)$ corresponding to eigenvalues with a 
negative real part. 
Then to construct the modified Evans function~\eqref{newdefevans} 
at the matching point $x_*$ we proceeded as follows (for both
the planar and wrinkled front cases for the autocatalytic problem below).
The full system is of size $n=4(2K+1)$.
We use the coordinate patch identified by $\ib_-=\{1,\ldots,m\}$
where $m=2(2K+1)$ for the interval $(-\infty,x_*]$. 
On this interval we solve the 
Riccati equation~\eqref{eq:riccati} for $\hat y_-$; with the chosen coordinate
patch, the coefficient matrices in the Riccati equation are
$a=O_{m}$, $b=I_{m}$, $c=A_3$ and $d=A_4$. The initial
data for $\hat y_-$ is simply the lower $(n-m)\times m$ block
of $Y_0^-(\lambda)$ postmultiplied by the inverse of the 
upper $m\times m$ block of $Y_0^-(\lambda)$. For the interval
$[x_*,+\infty)$ we use the coordinate patch identified by 
$\ib_+=\{m+1,\ldots,n\}$. We solve the Riccati equation~\eqref{eq:riccati} 
for $\hat y_+$; the coefficient matrices in this case are
$a=A_4$, $b=A_3$, $c=I_{m}$ and $d=O_{m}$. 
The initial data for $\hat y_+$ is the upper $m\times(n-m)$ block
of $Y_0^+(\lambda)$ postmultiplied by the inverse of the 
lower $(n-m)\times(n-m)$ block of $Y_0^+(\lambda)$.
Of course in practical calculations, our intervals of integration
are $[L_x^-,x_*]$ and $[x_*,L_x^+]$ for suitable choices of
$L_x^-$ and $L_x^+$. In particular, we assume the interval $[L_x^-,L_x^+]$ contains
the travelling front and extends far enough in both directions
that we can assume the front is sufficiently close to 
the far field homogeneous steady states. Further, 
the respective coordinate patches for the left and right intervals,
and matching point $x_*$, are chosen to avoid singularities in the Riccati flows 
in these intervals (in all our calculations for the autocatalysis problem we were
always able to find such a matching point when using these two patches).
Hence we evaluate the Evans function by computing the determinant~\eqref{newdefevans}. 
This information is used to find the spectrum of the travelling wave.

\subsection{The planar front}
\label{sec:planar}

We first consider planar travelling waves. This will not give us new
information but it is a good test for our procedure. The travelling
waves can be found by the freezing method, but also by the shooting
method explained in Balmforth, Craster and Malham~\cite{BCM} for the
one-dimensional case. 
The next step is to find the unstable subspace of~\eqref{eq:bigsys} at
$x=-\infty$. The travelling wave~$U_c$ is planar and thus independent
of~$y$, so the Fourier coefficient~$\hD_k$ vanishes for $k \ne 0$.
This implies that the differential equation~\eqref{eq:bigsys}
decouples in $2K+1$ systems, one for every~$k = -K, \ldots, K$, each
of the form
\begin{align*}
\partial_x\hU_k &= \hP_k, \\
\partial_x\hP_k &= \lambda B^{-1} \hU_k + \big(\tfrac{k}{\tL}\big)^2
\hU_k - B^{-1} \hD_0(x) \, \hU_k - c \, B^{-1} \hP_k.
\end{align*}
For the autocatalytic system~(\ref{cubicautoeqns}), this becomes
\begin{equation}
\label{eq:fourmatrix}
\partial_x \begin{pmatrix} \hU_k \\ \hP_k \end{pmatrix} = 
\begin{pmatrix} 
0 & 0 & 1 & 0 \\
0 & 0 & 0 & 1 \\
\tfrac\lambda\delta + \big(\tfrac{k}{\tL}\big)^2 + \tfrac1\delta u_2^2 
  & \tfrac2\delta u_1 u_2 & -\tfrac{c}\delta & 0 \\
-u^2_2 & \lambda + \big(\tfrac{k}{\tL}\big)^2 - 2u_1 u_2 & 0 & -c
\end{pmatrix}
\begin{pmatrix} \hU_k \\ \hP_k \end{pmatrix}.
\end{equation}
We have $(u_1,u_2) \to (0,1)$ in the limit $x\to-\infty$. The
coefficient matrix in~\eqref{eq:fourmatrix} has two unstable eigenvectors when
$\lambda$ is to the right of the continuous spectrum, and these
eigenvectors are
\begin{subequations}
\label{norm1}
\begin{align}
\left( 1, \tfrac1\nu, \mu, \tfrac\mu\nu \right)^{\mathrm{T}},
&\qquad\text{where}\qquad
\begin{cases}
\mu = -\frac{c}{2\delta} + \sqrt{\frac{c^2}{4\delta^2} +
    \big( \frac{k}{\tL} \big)^2 + \frac{\lambda+1}{\delta}}, \\
\nu = \big( \frac{k}{\tL} \big)^2 + \lambda - c\mu - \mu^2,
\end{cases}\\
\intertext{and}
(0, 1, 0, \mu)^{\mathrm{T}},
&\qquad\text{where}\qquad
\textstyle \mu = -\frac12 c + \sqrt{\frac14 c^2 + 
  \big( \frac{k}{\tL} \big)^2 + \lambda}.
\end{align}
\end{subequations}
These eigenvectors are all analytic functions of~$\lambda$ except at
those points where either $\nu$ or one of the expressions under the
square root sign vanishes. A small calculation shows that the latter
can only happen when $\lambda$ is on the negative real axis, so it
will not interfere with stability computations, which concern
eigenvalues with positive real part. In contrast, we may have $\nu=0$
for positive $\lambda$, but in the situations considered here $\nu$
vanishes only for values of $\lambda$ which are so large that they do
not concern us here. Alternatively, the singularities caused by
$\nu=0$ can easily be avoided by using $(\nu, 1, \mu\nu, \mu)^{\mathrm{T}}$ as
eigenvector instead of $(1, 1/\nu, \mu, \mu/\nu)^{\mathrm{T}}$.

We thus have two unstable directions for every wavenumber~$k$.
Together these form the $2(2K+1)$-dimensional unstable subspace
$Y_0^-(\lambda)$ at $x=-\infty$. Using this, as described in
Section~\ref{sec:method} above, we integrate the 
Riccati equation~\eqref{eq:riccati} from $x=L_x^-$ to $x=x_*$, 
where $L_x^-$ corresponds to a position sufficiently far behind the front. 
The integration is done with Matlab function
\verb|ode45|, which is based on the explicit fifth-order Runge--Kutta
method due to Dormand and Prince. In all experiments, we use an
absolute tolerance of~$10^{-8}$ and a relative tolerance
of~$10^{-6}$.

In the other limit, where $x\to+\infty$, we need the two stable
eigenvectors of the coefficient matrix in~\eqref{eq:fourmatrix}, which are
\begin{subequations}
\label{norm2}
\begin{align}
(1, 0, \mu, 0)^{\mathrm{T}},
&\qquad\text{where}\qquad
\textstyle \mu = -\frac{c}{2\delta} - \sqrt{\frac{c^2}{4\delta^2} +
    \big( \frac{k}{\tL} \big)^2 + \frac{\lambda}{\delta}},\\
\intertext{and}
(0, 1, 0, \mu)^{\mathrm{T}},
&\qquad\text{where}\qquad
\textstyle \mu = -\frac12 c - \sqrt{\frac14 c^2 + 
  \big( \frac{k}{\tL} \big)^2 + \lambda}.
\end{align}
\end{subequations}
Again, we get a $2(2K+1)$-dimensional subspace $Y_0^+(\lambda)$. 
We use this to generate the initial condition for the 
Riccati equation~\eqref{eq:riccati} we integrate backwards 
from $x=L_x^+$ to $x=x_*$---note the coefficients
are different to those for the left interval---see 
Section~\ref{sec:method} above. Here $L_x^+$ is a longitudinal position
sufficiently far ahead of the front. In our actual calculations
we used $L_x^\pm=\pm25$ and $x_*=0$.
Finally, the Evans function is computed using the modified form~\eqref{newdefevans}.

This concludes the explanation for the Riccati approach. The process
for the Drury--Oja approach due to Humpherys and Zumbrun~\cite{HZ} is
very similar. The difference is that it uses the Drury--Oja flow~\eqref{eq:DO} to
propagate the subspaces and~\eqref{eq:HZEvans} to evaluate the Evans
function.

\begin{figure}
\begin{center}
  \includegraphics[width=\linewidth]{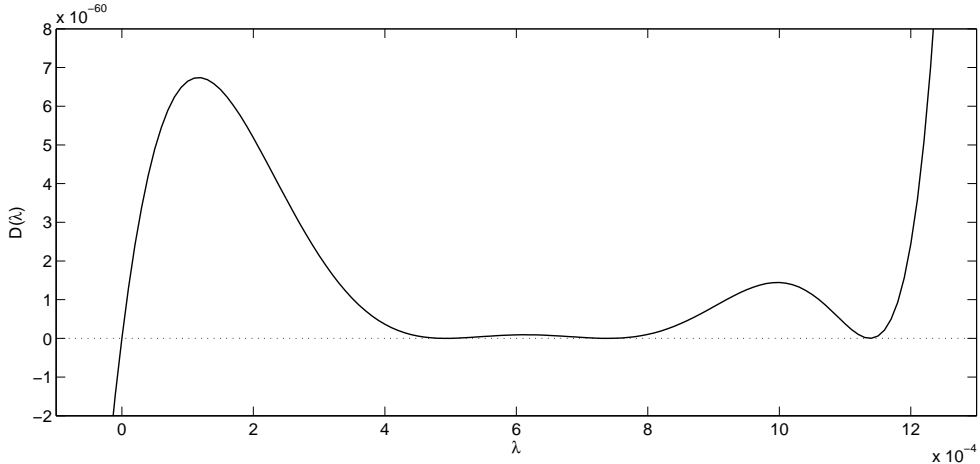} 
\end{center}
\caption{The Evans function along the real axis for the planar front
  at $\delta=3$. The Fourier cut-off is $K=24$ and the transversal
  length of the domain is $L=200$. The Evans function is computed
  using the Drury--Oja approach. A very similar plot can be produced
  using the Riccati approach.}
\label{evansp3}
\end{figure}

As an illustration, we show in Figure~\ref{evansp3} a graph of the
Evans function along the real axis for $\delta=3$, with the Fourier
cut-off at $K=24$. We see that the Evans function has a zero at
$\lambda = 0$ and double zeros around $0.0005$, $0.0007$ and $0.0011$.
The zero at $\lambda = 0$ corresponds to the eigenvalue at the origin
caused by translational symmetry. The other zeros correspond to double
eigenvalues. They have positive real part and thus we can conclude
that the planar wave is unstable for $\delta=3$---coinciding with
the long-established statements in Section~\ref{sec:auto}.

In fact, the zeros for the Evans function can be related to the
dispersion curve in Figure~\ref{dispersion}. As mentioned above, the
differential equation~\eqref{eq:bigsys} decouples into $2K+1$ subsystems
of the form~\eqref{eq:fourmatrix}, each corresponding to the linear
operator~$\cL$ given in~\eqref{eq:cL1d} for a particular wave number.
The basis vectors we chose for the unstable subspace are zero for all
but one of these subsystems. Since the subsystems are decoupled, the
solution at $x=x_*$ is also zero for all but one of these subsystems.
Thus, the matrix $\bigl(Y^-(x_*;\lambda)\,\, Y^+(x_*;\lambda)\bigr)$
at the matching point is (after reordering) a block-diagonal matrix
consisting of $2K+1$ four-by-four blocks, and its determinant is the
product of the determinants of the $4\times4$ blocks. However, every
determinant of a sub-block is the Evans function~$D_k$ of the
operator~\eqref{eq:cL1d} for a particular wave number~$k$. Thus, the
Evans function for the two-dimensional planar wave is the product of
the Evans functions for the one-dimensional wave with respect to
transverse perturbations:
\begin{equation}
\label{eq:factorization}
D(\lambda) = \prod_{k=-K}^K D_k(\lambda).
\end{equation}
We chose $L=200$, so $k=1$ corresponds to a wave number of
$\frac{2\pi}{200} \approx 0.0314$ in Figure~\ref{dispersion}. The plot
shows that the corresponding growth rate is approximately 0.0005, thus
$D_1(\lambda)$ has a zero around $\lambda = 0.0005$. The dispersion
relation is symmetric, so
$D_{-1}(\lambda)$ has a zero at the same value. This explains why the
two-dimensional Evans function has a double zero around $\lambda =
0.0005$. Similarly, the double zeros around $\lambda = 0.0011$ and
$\lambda = 0.0007$ correspond to $k=2$ and $k=3$, respectively.

As Figure~\ref{evansp3} shows, the Evans function is of the
order~$10^{-60}$ in the interval of interest. It is important
to emphasize that we are primarily interested in the zeros of the 
Evans function, which correspond to eigenvalues of the operator~$\mathcal{L}$ 
defined in~\eqref{eq:call}. The scale of the Evans function between the zeros,
which depends on the normalization of the eigenvectors~\eqref{norm1}
and~\eqref{norm2}, is of lesser relevance in this context. 
However for completeness, we address this issue in Section~\ref{sec:convergence}.

\subsection{The wrinkled front}
\begin{figure}
\begin{center}
\includegraphics[width=0.48\linewidth]{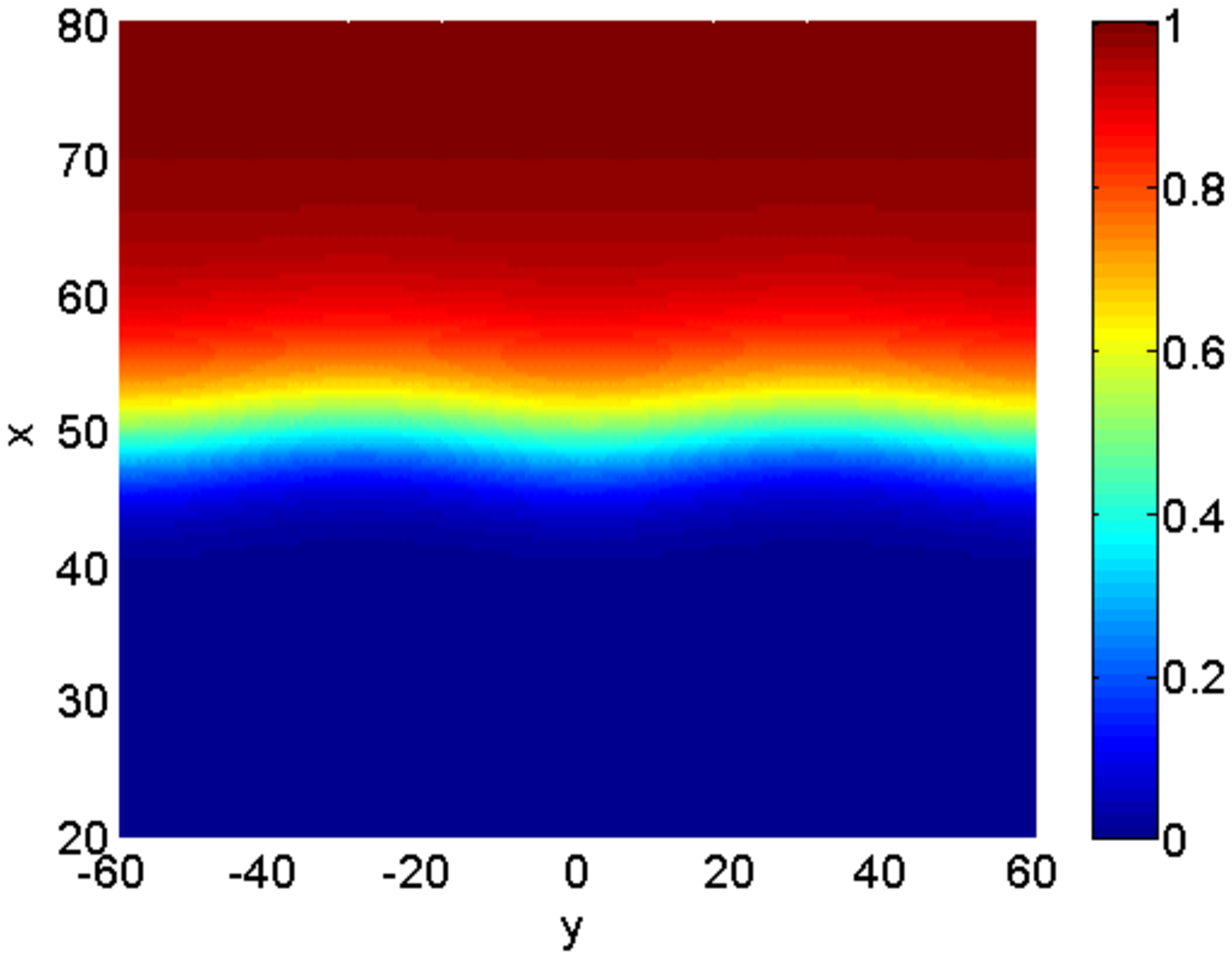}
\includegraphics[width=0.48\linewidth]{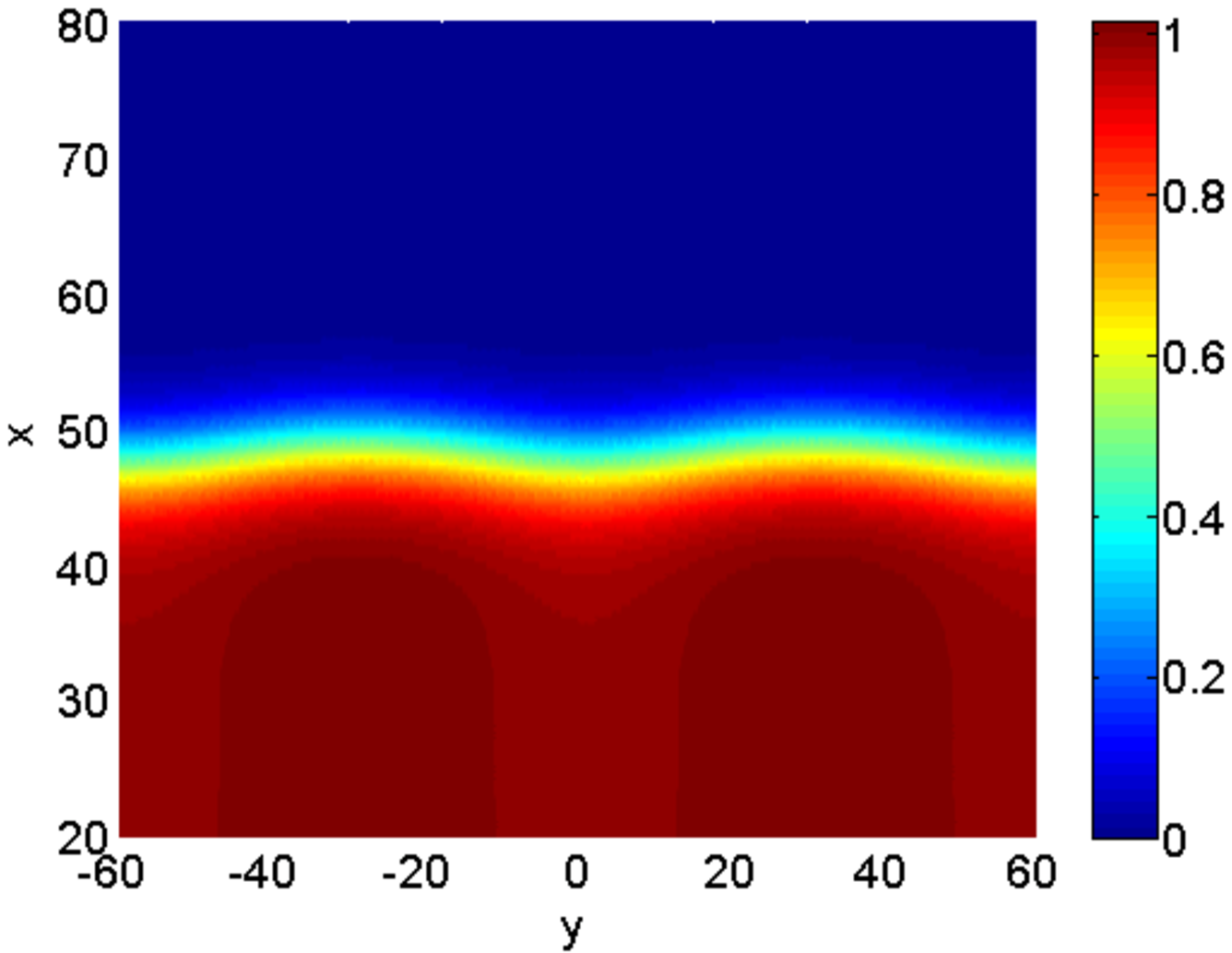}
\end{center}
\caption{The wrinkled front for $\delta=3$. The left panel shows the
  $u$~component and the right panel shows the $v$~component.}
\label{front}
\end{figure}

Having established that our method works for planar fronts, we now turn to the
wrinkled front. The most immediate problem is that of computing reliable
steady wrinkled fronts. The procedure we used to do this is outlined as
follows:
\begin{itemize}
\item The partial differential algebraic system~\eqref{eq:pdae} 
is solved on the computational domain $[-150,150] \times [-60,60]$ with a
second-order finite element method on a grid of right triangles
formed by dividing a grid of $300 \times 240$ rectangles along their
diagonals. We used the Finite Element package Comsol
Multiphysics\texttrademark~\cite{Comsol} which includes a version of
the DAE solver~\texttt{daspk}.
\item To obtain reasonable starting waveform profiles for
the partial differential algebraic system~\eqref{eq:pdae},
we constructed the one-dimensional waveform, which we then swept
out uniformly in the transverse direction.
\item For travelling waves that appeared to be stable, we use the
freezing method described in Section~\ref{sec:comp}
to obtain good starting waveform profiles to initiate the use of Newton's 
method in Comsol. 
\item For travelling waves that appeared to be unstable, we use simple 
parameter continuation (also using Comsol) starting with 
a (parametrically) nearby stable wave.
\end{itemize}
Figure \ref{front} shows the front which results for $\delta=3$.

The wrinkled front varies along the transverse $y$-direction.
Therefore, the Fourier coefficients~$\hD_k$ do not vanish for $k \ne
0$ and the differential equation~\eqref{eq:bigsys2} does not decouple.
In the limits $x \to \pm\infty$, we have restricted ourselves to
wrinkled front profiles that approach a $y$-independent state, 
which is the same as for the planar front. 
Hence, the computation of the stable and unstable
subspaces for the planar front, which we reported above, remains valid
for the wrinkled front.

\begin{figure}
\includegraphics[width=.98\linewidth]{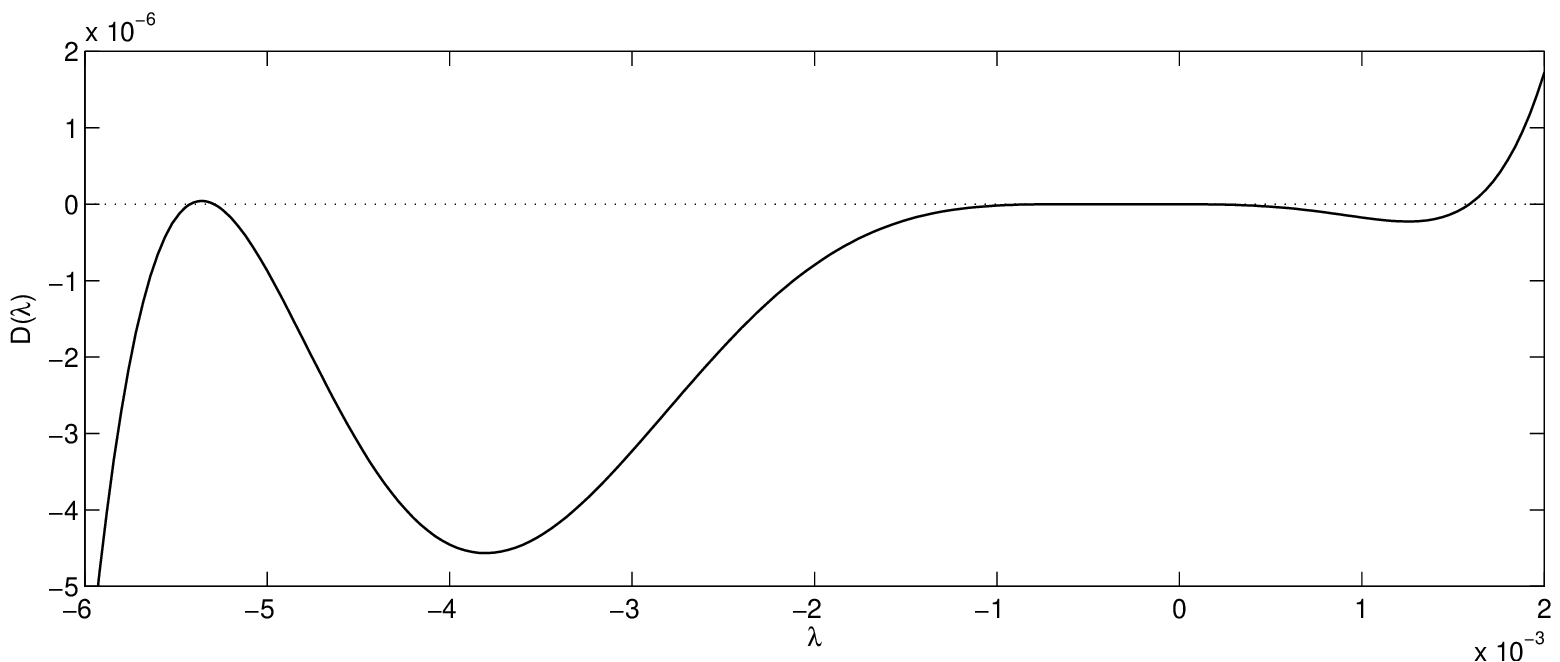}\\ 
\makebox[\linewidth]{
  \includegraphics[width=0.31\linewidth]{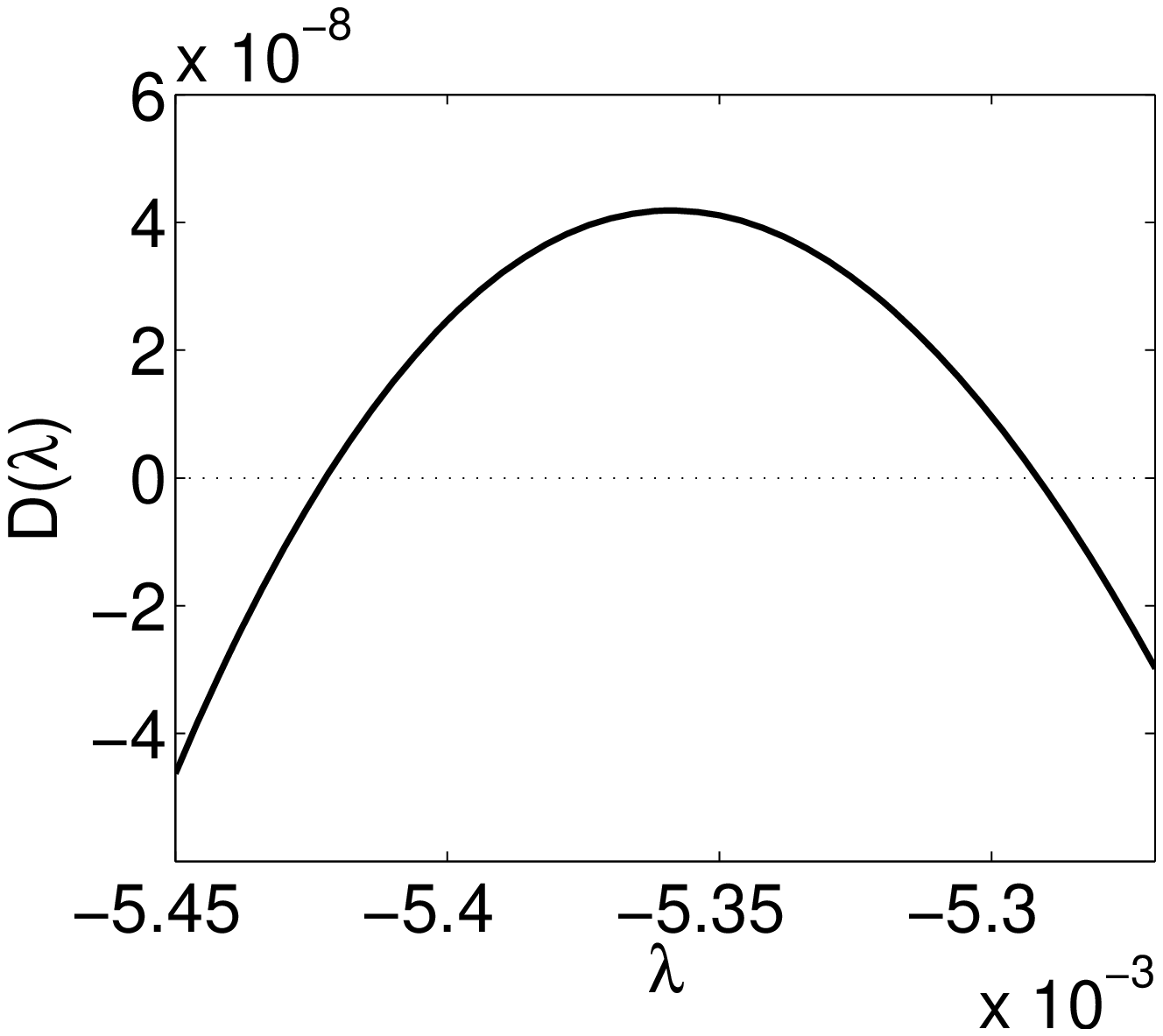}\hfill
  \includegraphics[width=0.31\linewidth]{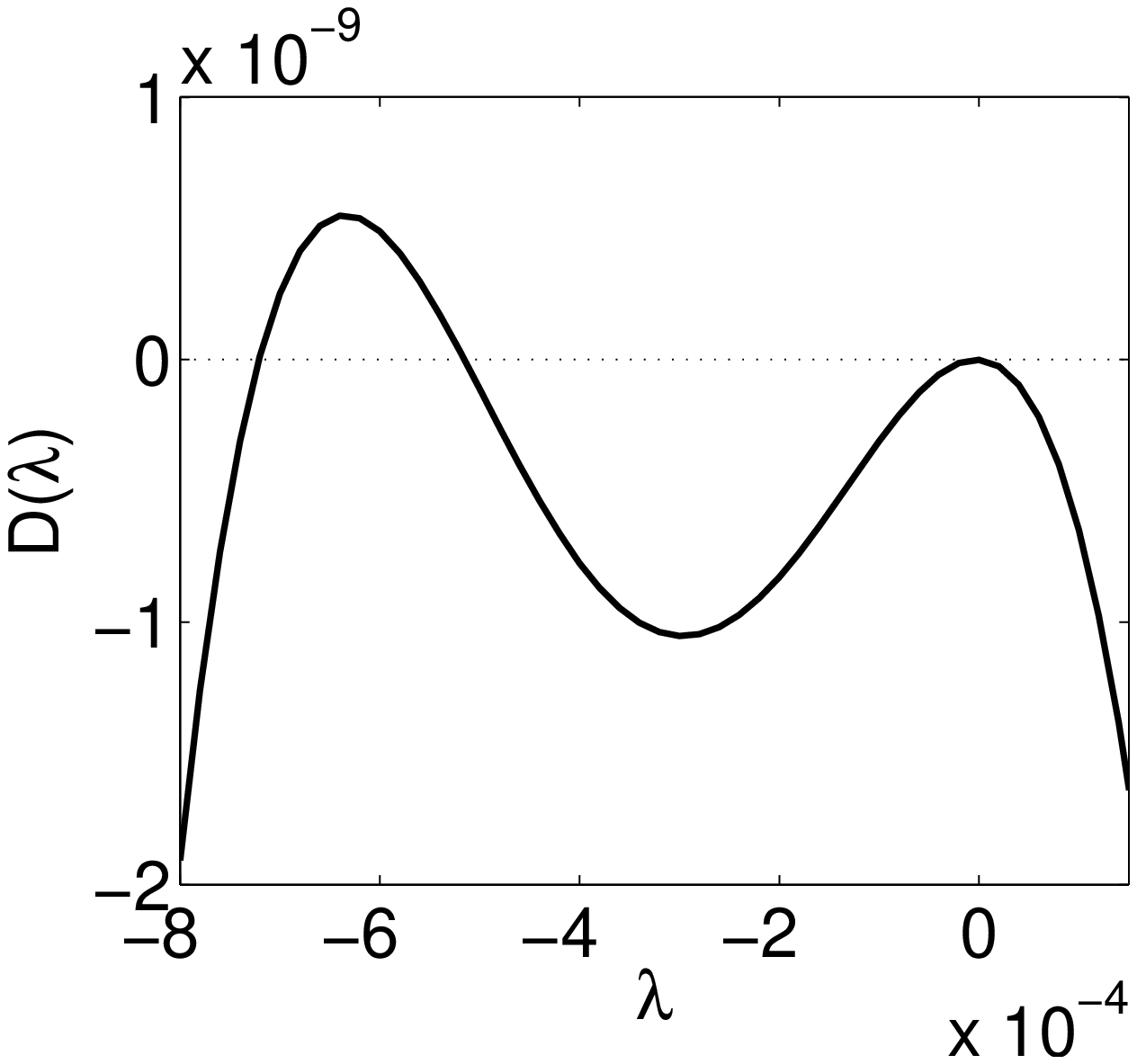}\hfill
  \includegraphics[width=0.31\linewidth]{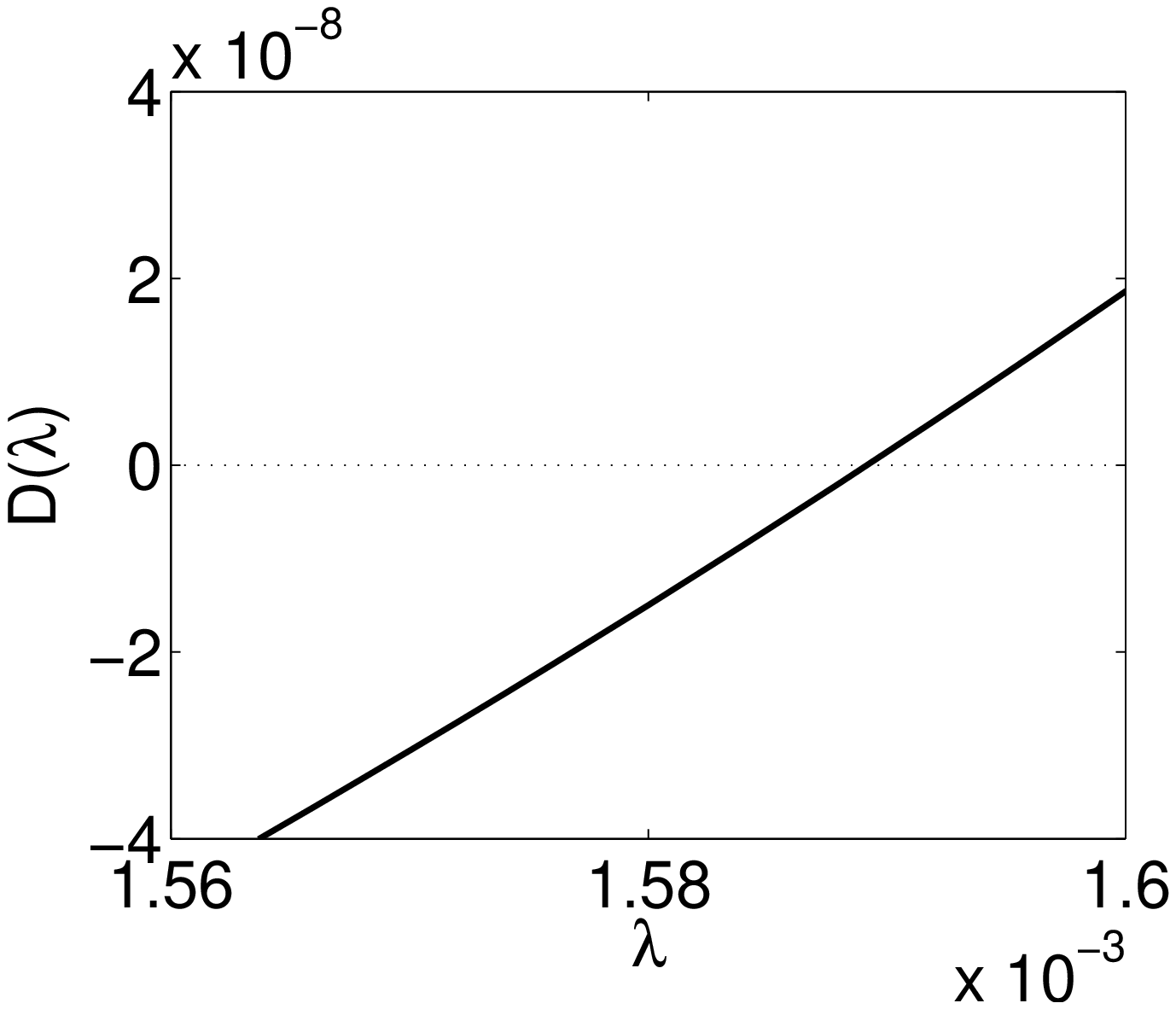}
}
\caption{The Evans function along the real axis for the wrinkled front
  at $\delta=3$, computed using the Riccati approach. 
  The three plots in the bottom row zoom in on the
  zeros of the Evans function.}
\label{evansw3ric}
\end{figure}

\begin{figure}
\includegraphics[width=.98\linewidth]{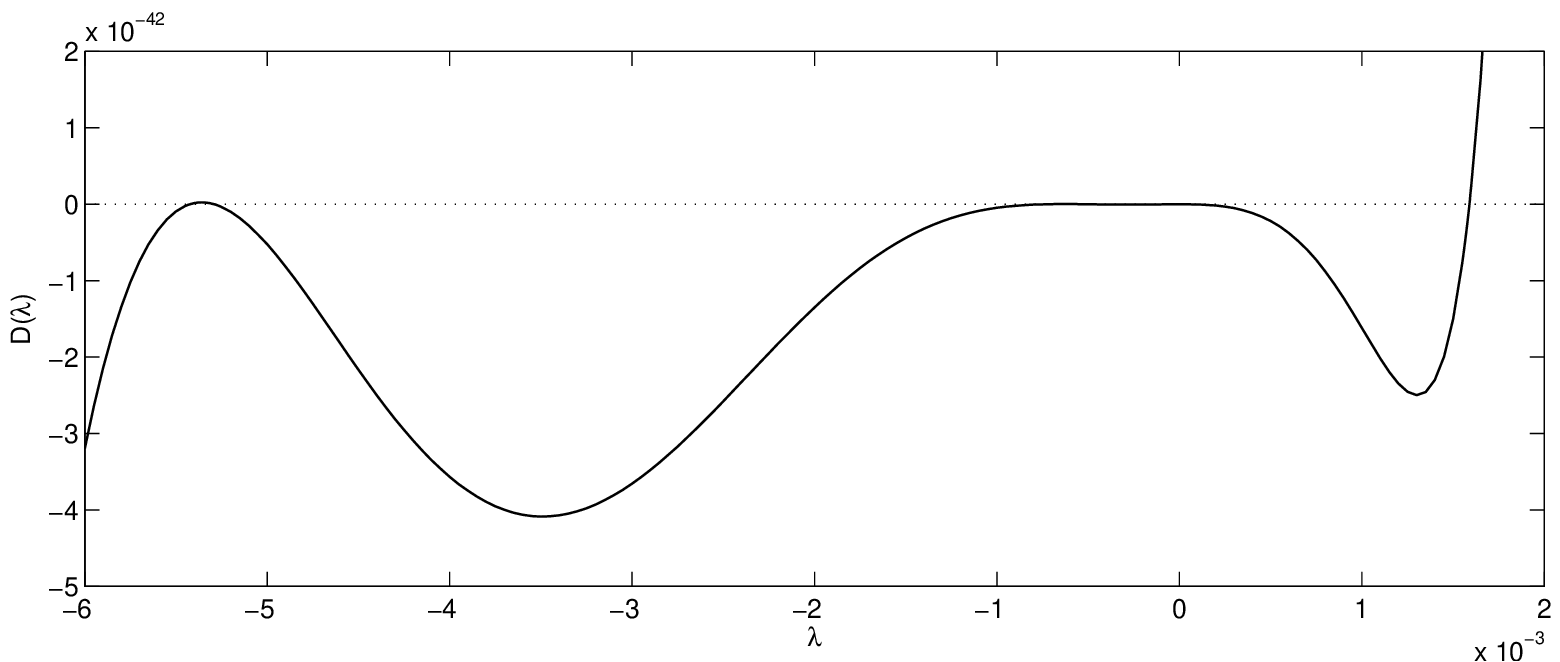}\\ 
\makebox[\linewidth]{
  \includegraphics[width=0.31\linewidth]{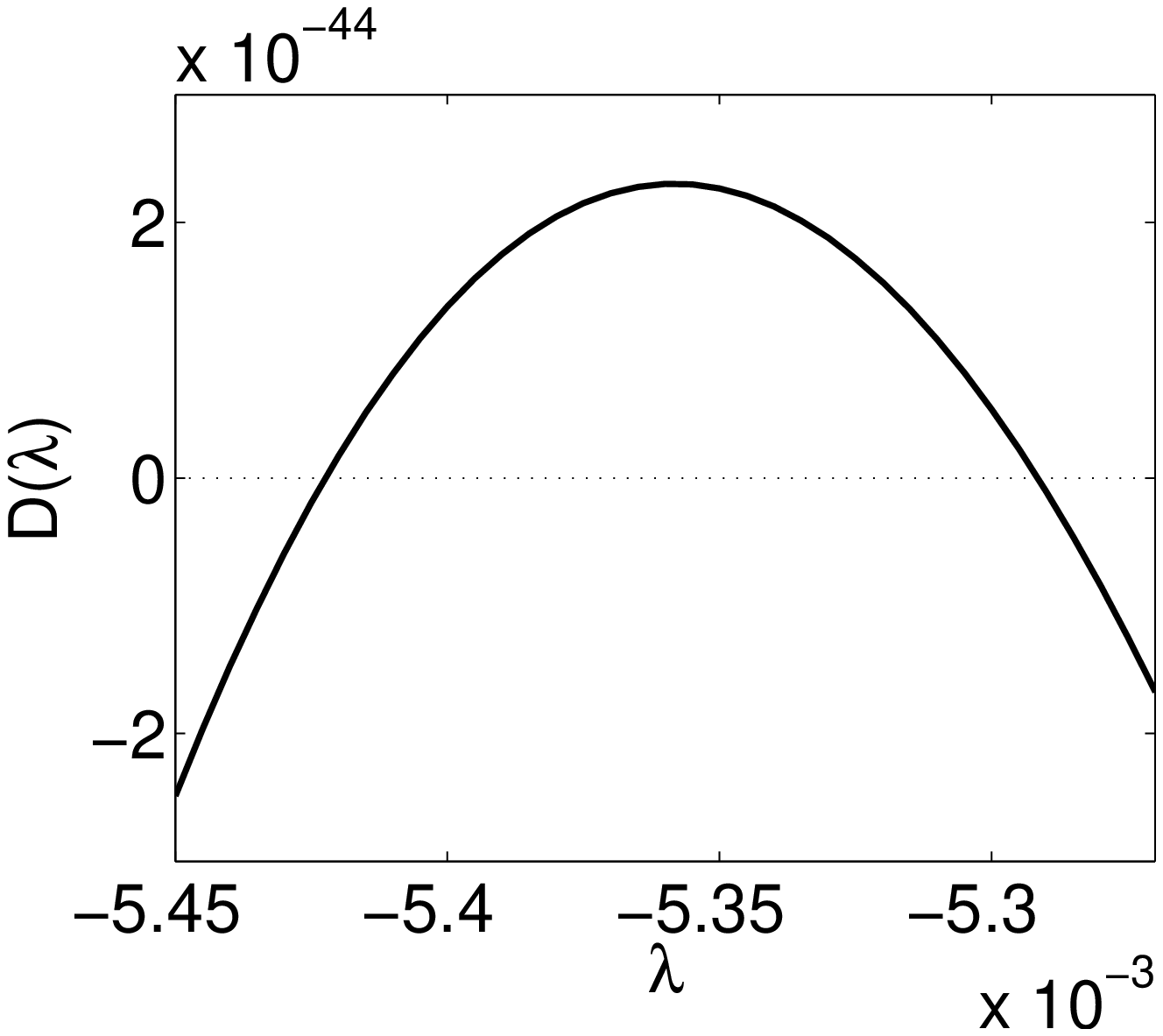}\hfill
  \includegraphics[width=0.31\linewidth]{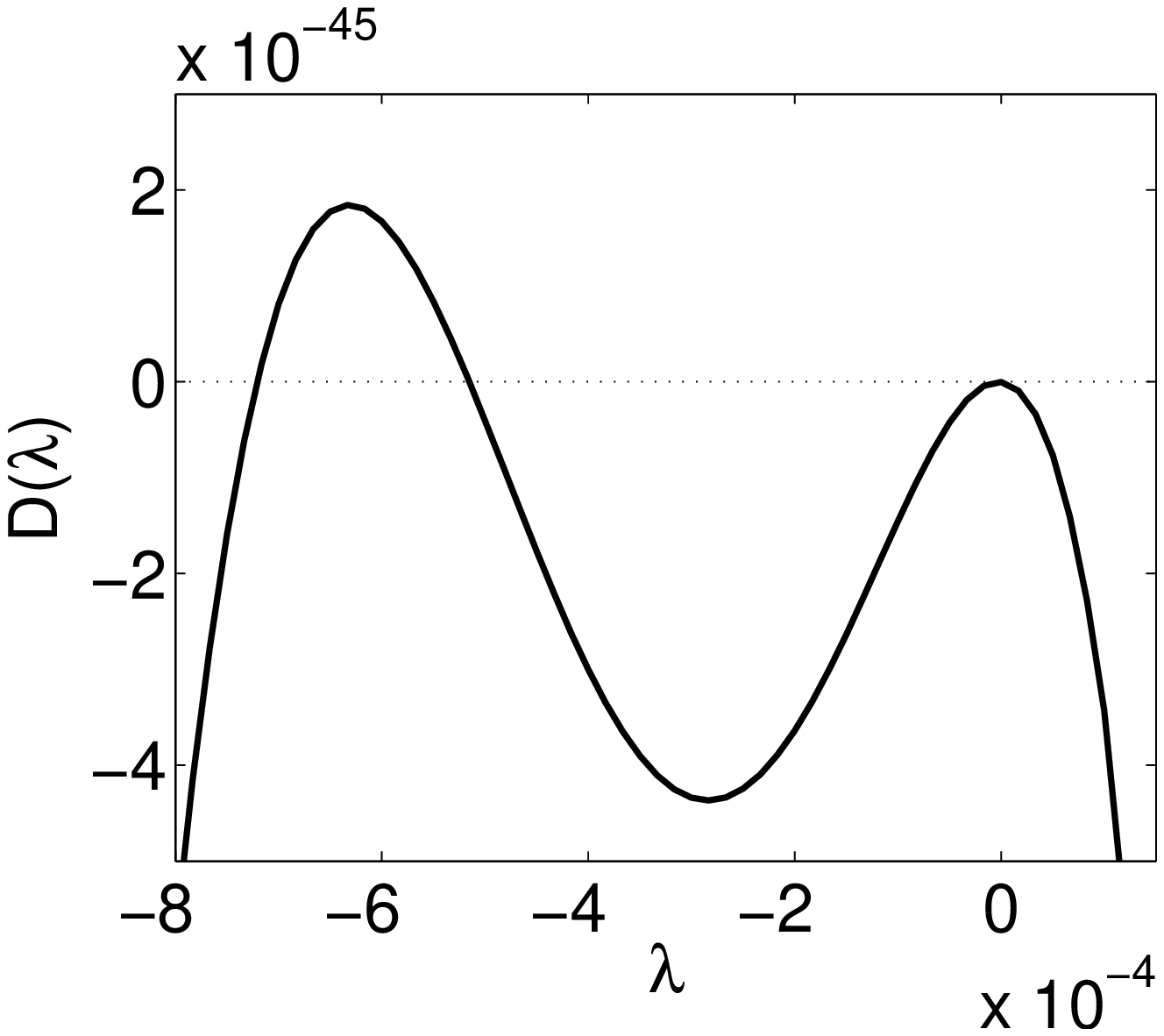}\hfill
  \includegraphics[width=0.31\linewidth]{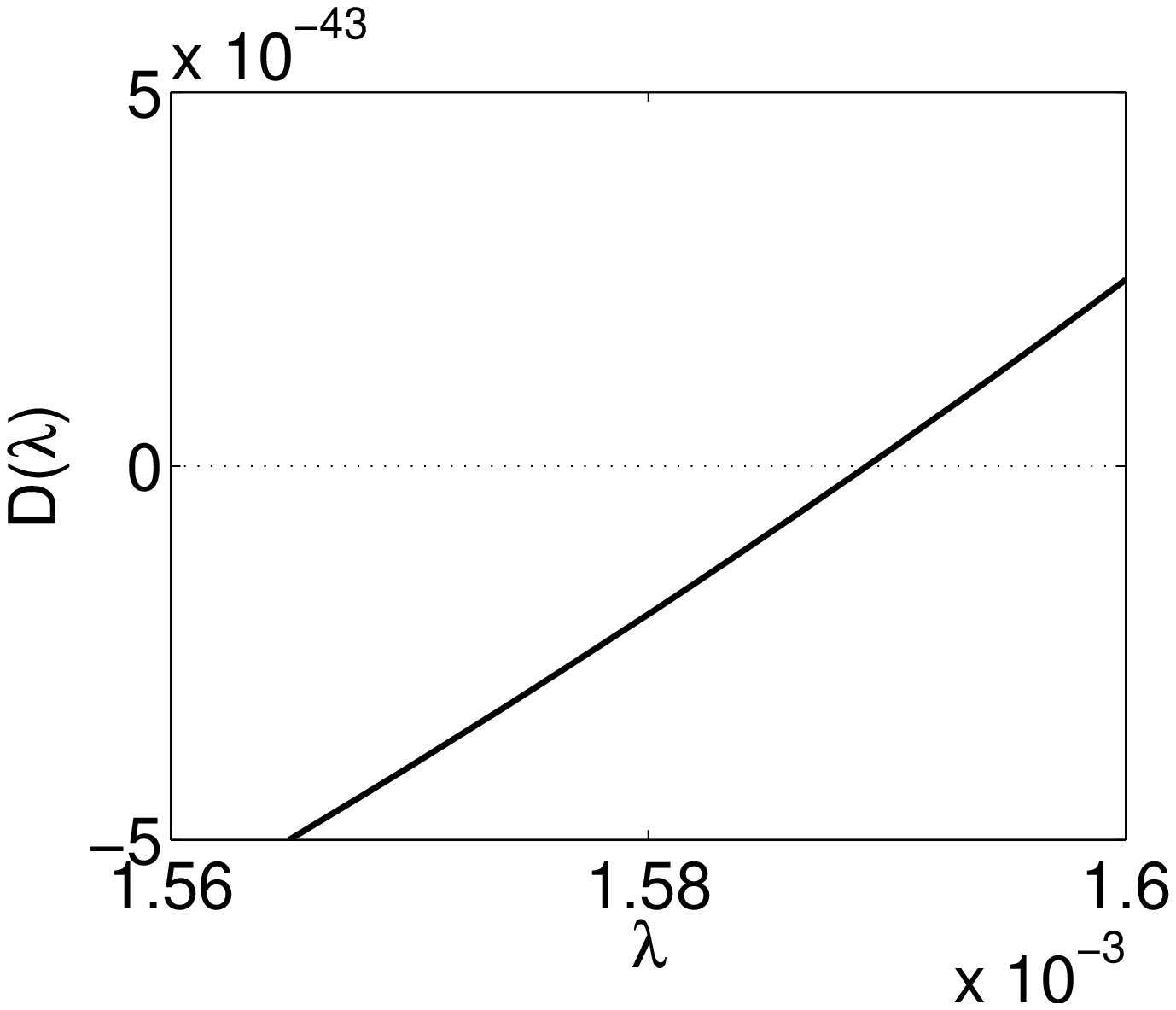}
}
\caption{Same as Figure~\ref{evansw3ric} but with the Drury--Oja
  approach.} 
\label{evansw3qr}
\end{figure}

The computation now proceeds as with the planar front. We solve the
Riccati differential equation~\eqref{eq:riccati} in the left and
right intervals (as described in Section~\ref{sec:method} above)
for the Riccati approach, and~\eqref{eq:DO} for the Drury--Oja approach. 
The front~$U_c$ is determined as explained above, with interpolation being
used for those points that fall between the grid points of the finite
element solution. We use cubic interpolation, because linear
interpolation is not sufficiently accurate. Finally, the Evans
function is computed using either the modified form~\eqref{newdefevans}
or Humpherys and Zumbrun form~\eqref{eq:HZEvans}.

Figure~\ref{evansw3ric} shows a plot of the Evans function along the
real axis for $\delta=3$, computed using the Riccati approach.  The
result of the Drury--Oja approach is shown in Figure~\ref{evansw3qr}.
For both approaches, we used projection on $K=9$ modes and we solved
the Riccati differential equations~\eqref{eq:riccati} (with appropriate coefficients
for the left and right intervals),
and the Drury--Oja equations~\eqref{eq:DO}, 
by the Matlab ODE-solver \verb|ode45| with absolute
and relative tolerances of $10^{-8}$ and $10^{-6}$. The domain is
$[-150,150]\times[-60,60]$ and as matching point we have chosen
$x_*=50$ which is roughly the front location (see Figure~\ref{front}). 
Both approaches produce Evans functions with zeros at (approximately)
the same $\lambda$-values: there is one eigenvalue with a positive
real part around $\lambda=0.0016$, a double zero around $0$ (due to
translation invariance), and pairs of eigenvalues around~$-0.0006$
and~$-0.0053$. In particular, the eigenvalue around~$\lambda=0.0016$
shows that the wrinkled front is unstable for $\delta = 3$.

However, the Riccati and Drury--Oja approaches do differ in some
respects. The Riccati approach involves integrating a system of half
the dimension as compared to the Drury--Oja approach, so we can 
expect it to be faster, especially for large values of $K$ 
(in comparable experiments for large $K$, it was typically two to three
times faster).
See Ledoux, Malham and Th\"ummler~\cite{LMT} where some detailed
accuracy versus cost comparisons can be found. 
We also note that the Drury--Oja approach again produced 
very small values for the Evans function, 
which as explained in Section~\ref{sec:convergence}, is
not as alarming as it first seems. The Riccati approach, 
which uses a different scaling of the Evans function, yields only
moderately small values, though again see Section~\ref{sec:convergence}.

We also used the angle between the unstable and stable subspace
$\theta(\lambda;x_\ast)$ to determine the location of the eigenvalues.
Though this scales better than the Evans determinant, it is
non-negative by definition and zeros occur as touch-downs in
the complex spectral parameter plane. Actual zeros of the function are
thus harder to establish compared to sign changes in the real and
imaginary parts of the Evans determinant. Further,
$\theta(\lambda;x_\ast)$ is not analytic in~$\lambda$.

\begin{table}
\begin{tabular}{@{}c|c@{ \ }r@{ \ }c@{ \ }c@{ \ }c@{ \ }c@{ \ }c@{}}
\hline
$K$ & \multicolumn{7}{c}{Eigenvalues (Evans function)} \\ \hline
3 & $0.001609$ & $-0.000026$ & $-0.000781$ & $-0.001296$ & $-0.000670$ & $-0.006078$ & $-0.006177$ \\
4 & $0.001609$ &  $0.000002$ & $-0.000001$ & $-0.000519$ & $-0.000670$ & $-0.006078$ & $-0.006177$ \\
5 & $0.001589$ &  $0.000002$ & $-0.000001$ & $-0.000519$ & $-0.000720$ & $-0.005295$ & $-0.005426$ \\
6 & $0.001589$ & $-0.000002$ & $-0.000003$ & $-0.000515$ & $-0.000720$ & $-0.005295$ & $-0.005426$ \\
7 & $0.001589$ & $-0.000002$ & $-0.000003$ & $-0.000515$ & $-0.000721$ & $-0.005292$ & $-0.005422$ \\
8 & $0.001589$ & $-0.000002$ & $-0.000003$ & $-0.000515$ & $-0.000721$ & $-0.005292$ & $-0.005422$ \\
9 & $0.001589$ & $-0.000002$ & $-0.000003$ & $-0.000515$ & $-0.000721$ & $-0.005292$ & $-0.005422$ \\[-1ex]
$\vdots$ & \multicolumn{7}{c}{$\vdots$} \\ 
24& $0.001589$ & $-0.000002$ & $-0.000003$ & $-0.000515$ & $-0.000721$ & $-0.005292$ & $-0.005422$ \\
\hline\hline
 & \multicolumn{7}{c}{Eigenvalues (ARPACK)} \\ \hline
 &$0.001592$&$0.000000$&$0.000000$ & $-0.000514$ & $-0.000719$ & $-0.005289$ & $-0.005420$ \\
\hline
\end{tabular}
\medskip
\caption{Zeros of the Evans function for the wrinkled front at
  $\delta=3$ for different values of the wave number cut-off~$K$, as
  computed with the Riccati and Drury--Oja approaches. The last row
  shows the eigenvalues computed by ARPACK.}
\label{evanstable}
\end{table}

More accurate values for the zeros of the Evans function can be found
by using a standard root-finding method. This yields the values listed
in Table~\ref{evanstable}. The results of the Riccati and Drury--Oja
approaches agree (to the precision shown) in all cases. The table also
shows results for different values of the wave number cut-off~$K$. The
eigenvalues converge quickly when $K$ increases, and projection on $K
= 7$ modes is sufficient to get six digits of accuracy, at least in
this example. The largest computation we performed is with $K = 24$,
in which case the one-dimensional equation has order~$196$.

We also computed the eigenvalues by a projection method to check these
results. The Comsol package, which we used to compute the travelling
front, has an eigenvalue solver based on the Arnoldi method as
implemented in ARPACK~\cite{LSY}. This solver uses a shift~$\sigma$ to
allow the user to choose on which part of the spectrum to concentrate. 
We found that the choice of this shift is not straightforward. 
It is unwise to use a shift which is exactly an eigenvalue, 
so in particular, the default value of $\sigma=0$ should not be used. 
We used a shift of $\sigma=1$ or $\sigma=3$ which proved to be safe. 
The resulting eigenvalues are also listed in Table~\ref{evanstable}. 
They are in close agreement with the eigenvalues found using the Evans
function, giving additional confidence in our computation.

\begin{figure}
\begin{center}
\includegraphics[width=0.8\linewidth]{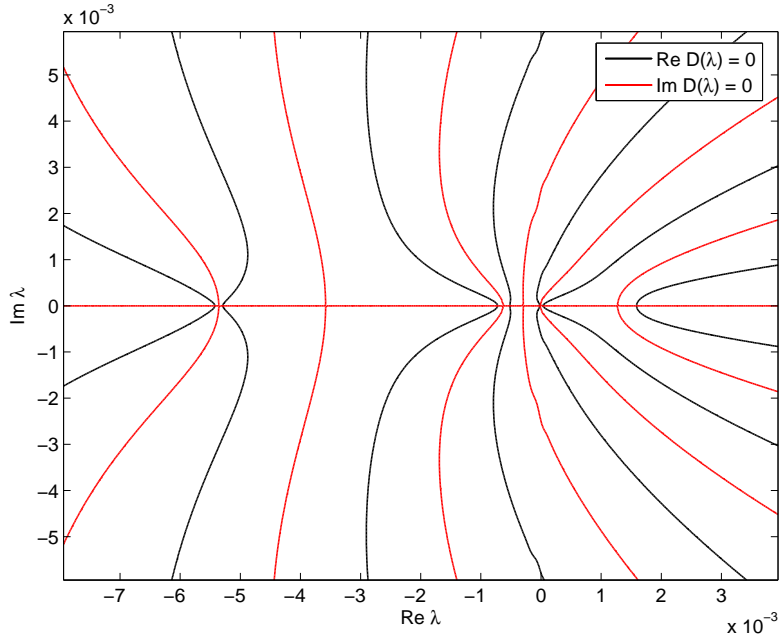}
\end{center}
\caption{The locus where the real and imaginary parts of Evans
  function vanish for the wrinkled front at $\delta=3$. This plot uses
  a $100\times100$ grid in the depicted region of the complex plane.
  The Evans function was computed on the 5000 grid points in the top
  half and extended by symmetry to the bottom half.}
\label{evansw3cp}
\end{figure}

Figure~\ref{evansw3cp} shows the contours where the real and imaginary
parts are zero; eigenvalues correspond with the intersections of these
contours. All eigenvalues are real. Interestingly, the
contours are well separated on the boundary of the depicted region of
the complex plane even though the eigenvalues are packed closely
together. This suggests that it may be better to study the Evans
function in a contour surrounding the eigenvalues instead of
concentrating on the eigenvalues themselves---see Section~\ref{sec:globalsearch}.

\begin{figure}
\begin{center}
\includegraphics[width=0.48\linewidth]{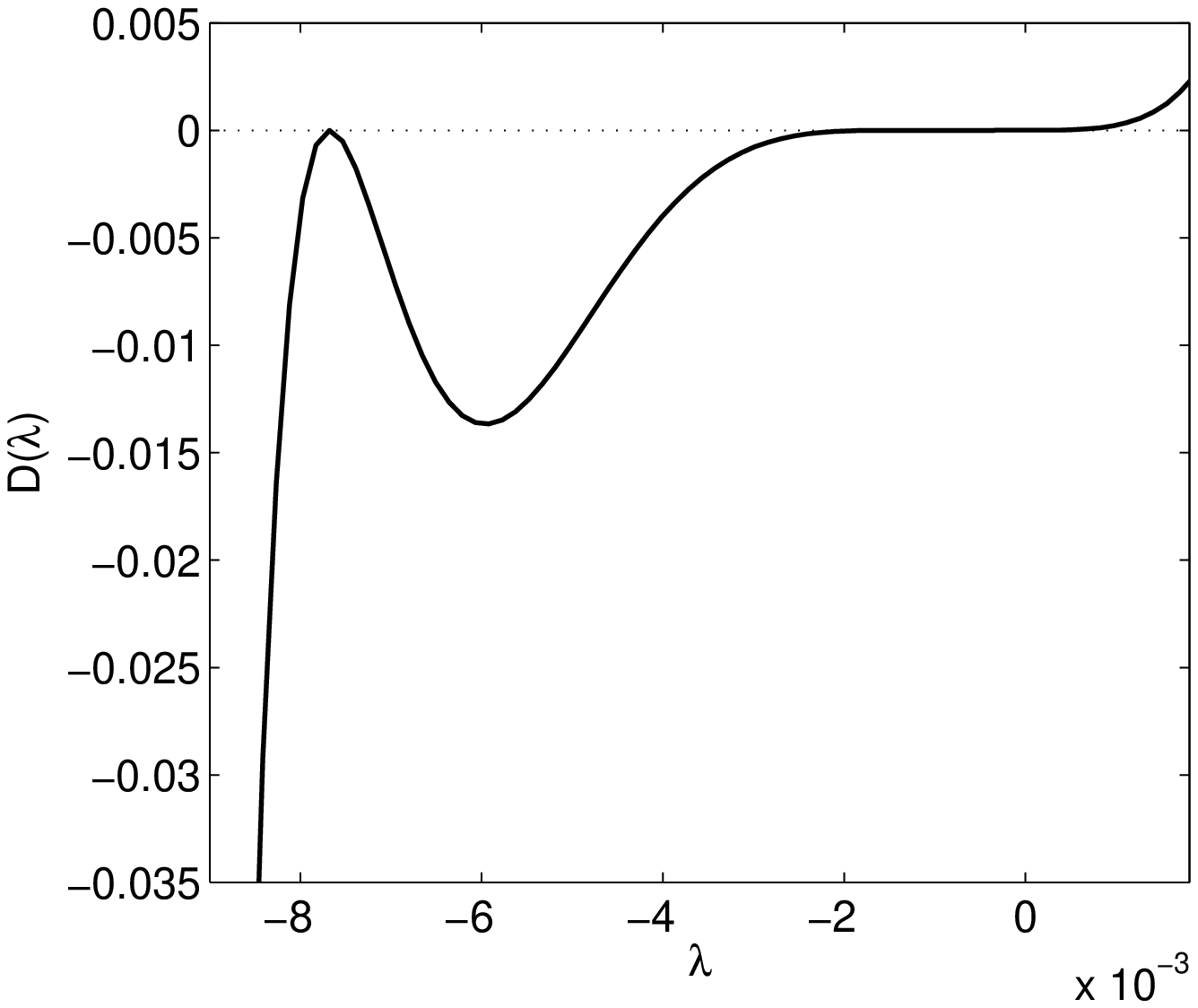} 
\includegraphics[width=0.48\linewidth]{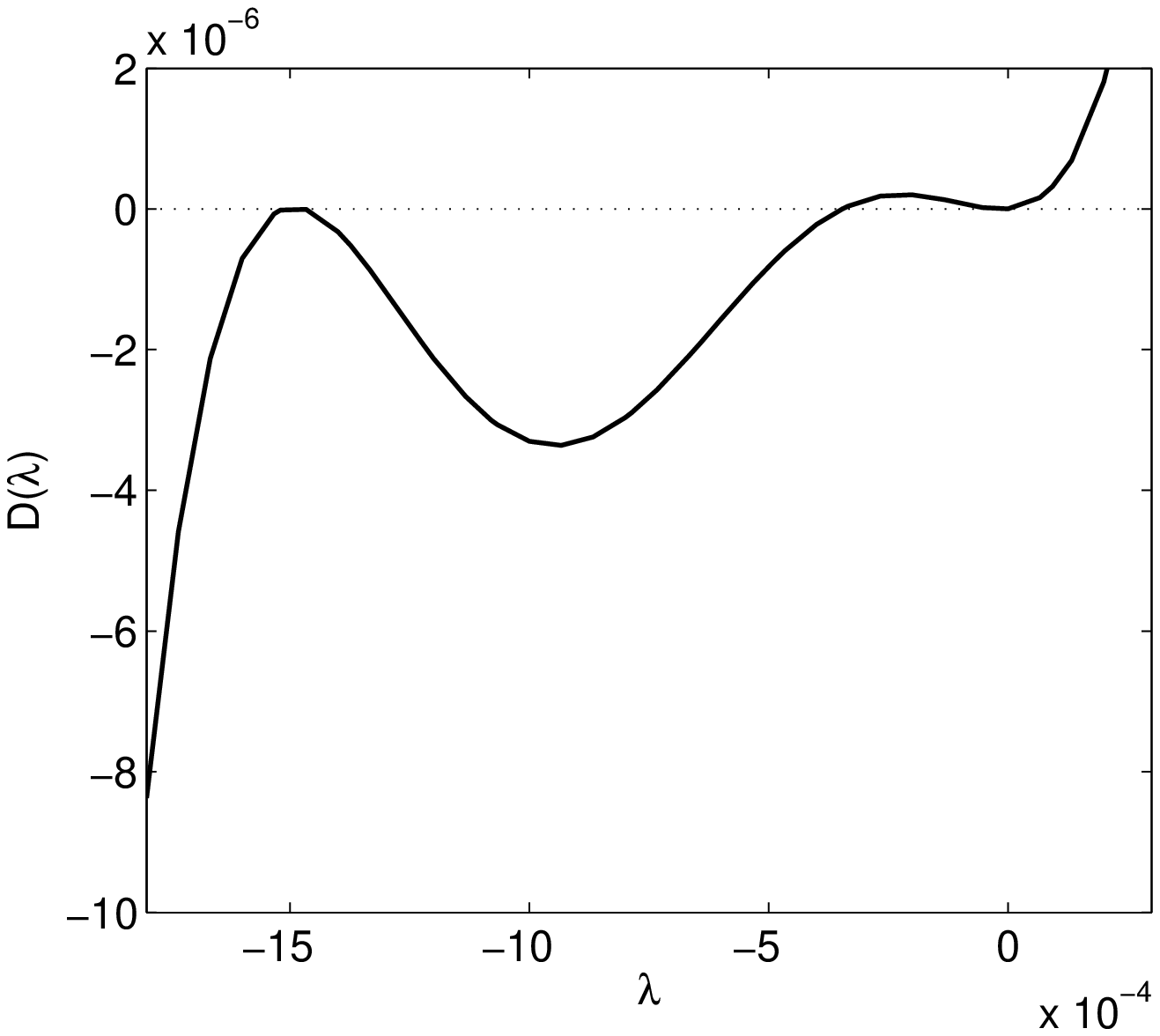}
\end{center}
\caption{The Evans function along the real axis for the wrinkled front
  at $\delta=2.5$, computed using the Riccati approach. 
  The right panel zooms in on the flat plateau around
  the origin in the left panel.}
\label{evansw2.5ric}
\end{figure}

\begin{figure}
\begin{center}
\includegraphics[width=0.48\linewidth]{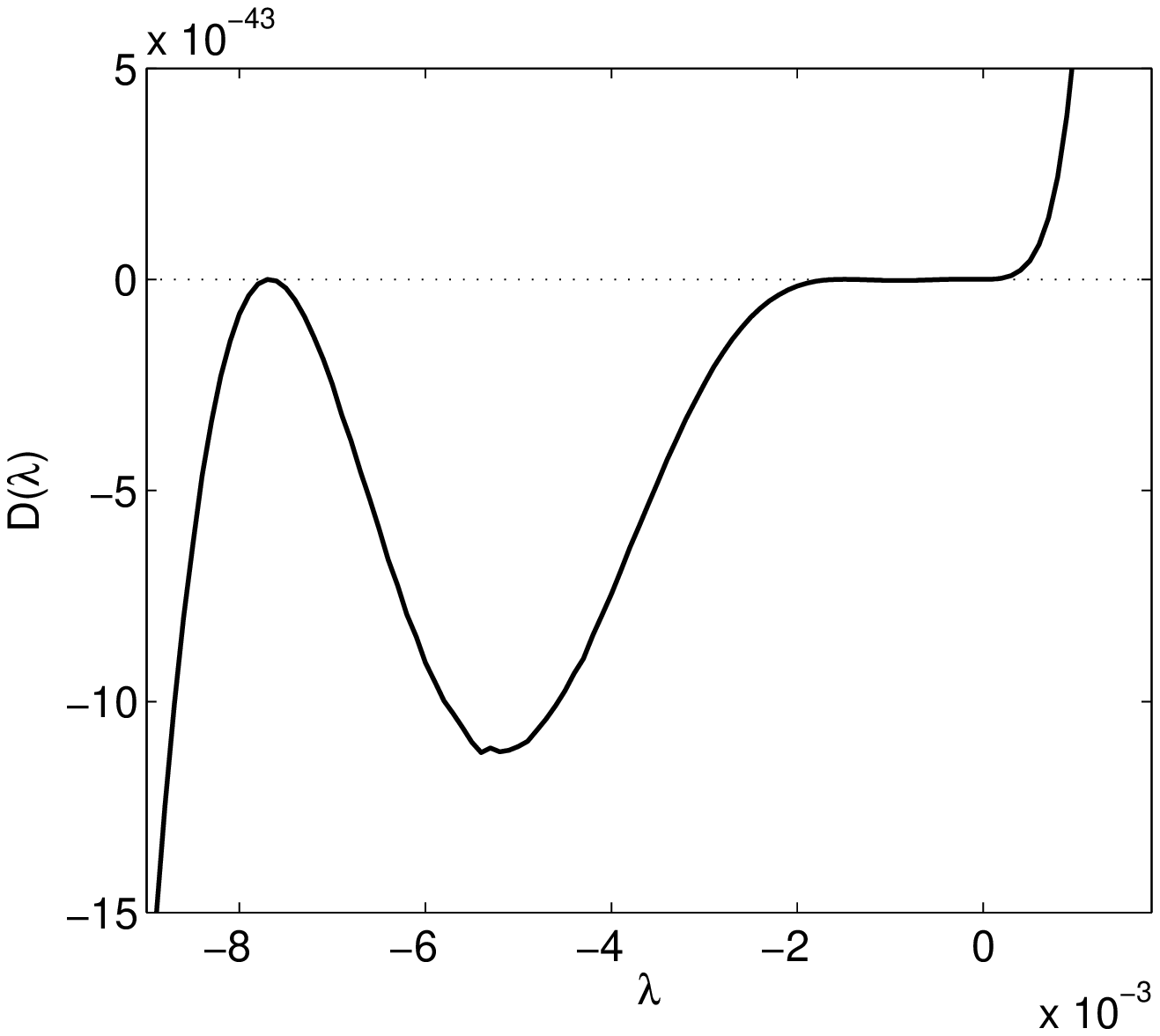} 
\includegraphics[width=0.48\linewidth]{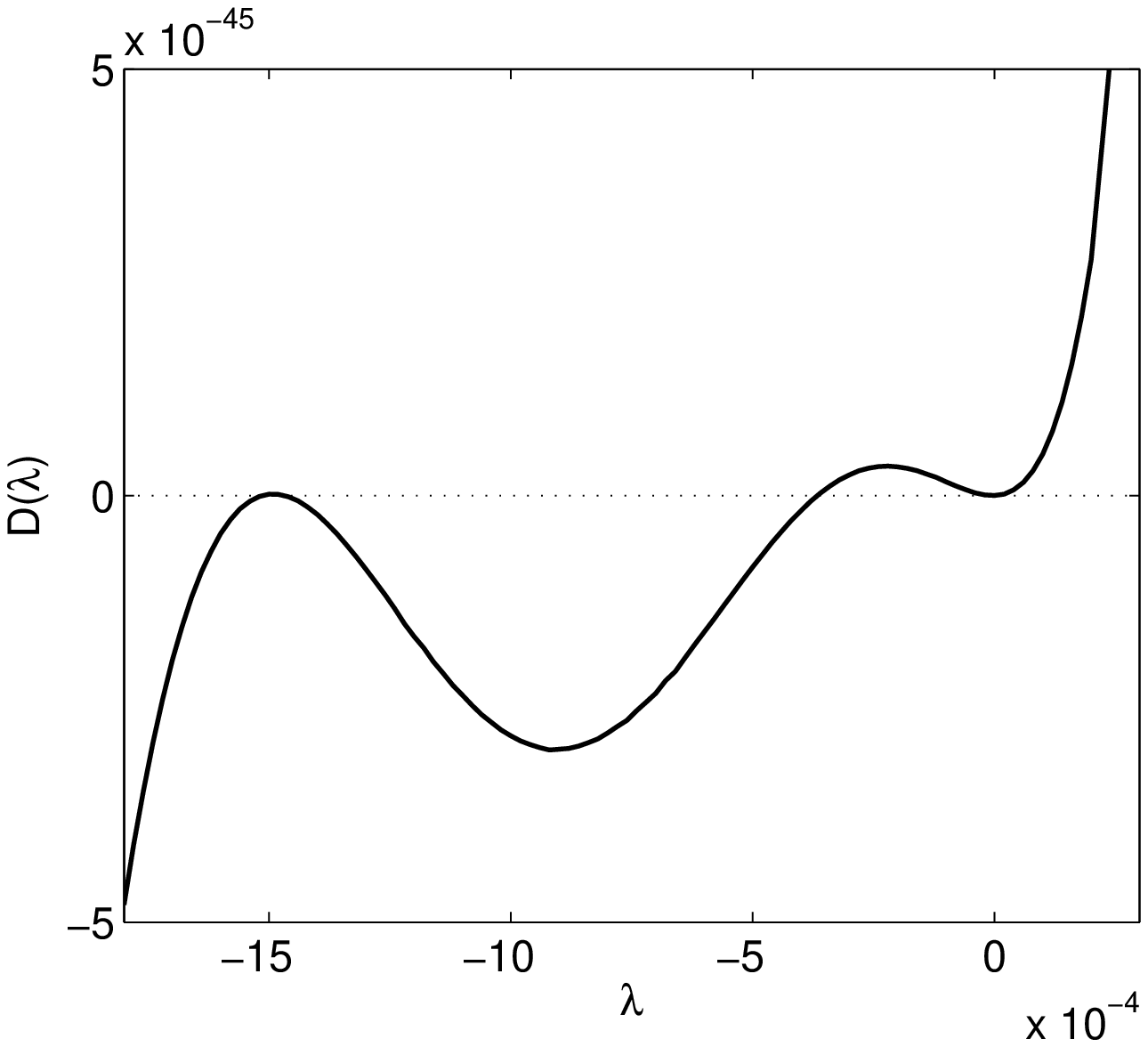}
\end{center}
\caption{Same as Figure~\ref{evansw2.5ric} but with the Drury--Oja
  approach.}
\label{evansw2.5qr}
\end{figure}

Figures~\ref{evansw2.5ric} and~\ref{evansw2.5qr} show the Evans
function for the case $\delta=2.5$. The Evans function no longer has a
zero with positive real part, meaning that there are no unstable
eigenvalues. Hence for this value of $\delta$ the wrinkled front
is in fact stable. We have thus established the following overall
stability scenario for travelling wave solutions of the cubic 
autocatalysis problem. When $\delta$ is less than 
$\delta_{\mathrm{cr}}\approx2.3$ planar
travelling waves are stable to transverse perturbations. 
For $\delta$ beyond this instability there exist wrinkled
fronts that are stable for $\delta=2.5$. However for $\delta=3$
these wrinkled fronts themselves become unstable. These stability
results may depend on the transverse domain size (which we fixed at
$L=120$); see Horv\'ath, Petrov, Scott and Showalter~\cite{HPSS}.

\subsection{Global eigenvalue search methods}\label{sec:globalsearch}
If one is only interested in the stability of the underlying front
then several methods exist to detect eigenvalues in the right-half
complex spectral parameter plane---which also require relatively
few Evans function evaluations. One important method uses the 
argument principle. We have the following sectorial estimate for
the spectrum of $\cL$: if $\lambda\in\sigma(\cL)$ then
\begin{subequations}
\label{eq:sectorial}
\begin{align}
\mathrm{Re}(\lambda)&
\leq\|\mathrm{D}F(U_c)\|_{L^\infty(\mathbb R\times\mathbb T)},
\label{eq:sectorial1} \\
\mathrm{Re}(\lambda)+|\mathrm{Im}(\lambda)|&\leq \frac{c^2}{4\kappa}
+2\,\|\mathrm{D}F(U_c)\|_{L^\infty(\mathbb R\times\mathbb T)},
\label{eq:sectorial2}
\end{align}
\end{subequations}
where $\kappa\equiv\min\{B_{11},B_{22}\}$---we provide a proof
in Appendix~\ref{app:sectorial}. Hence we can restrict our
search for unstable eigenvalues to the region of the complex
$\lambda$-plane bounded by these inequalities.

\begin{figure}
\begin{center}
\parbox[c]{0.24\linewidth}{\includegraphics[width=\linewidth]{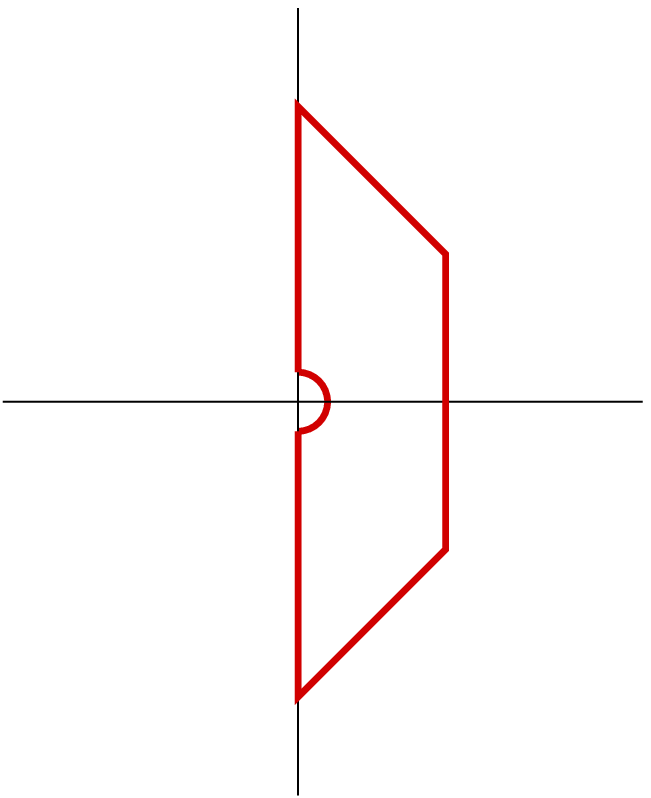}}%
\hspace{0.1\linewidth}
\parbox[c]{0.64\linewidth}{\includegraphics[width=\linewidth]{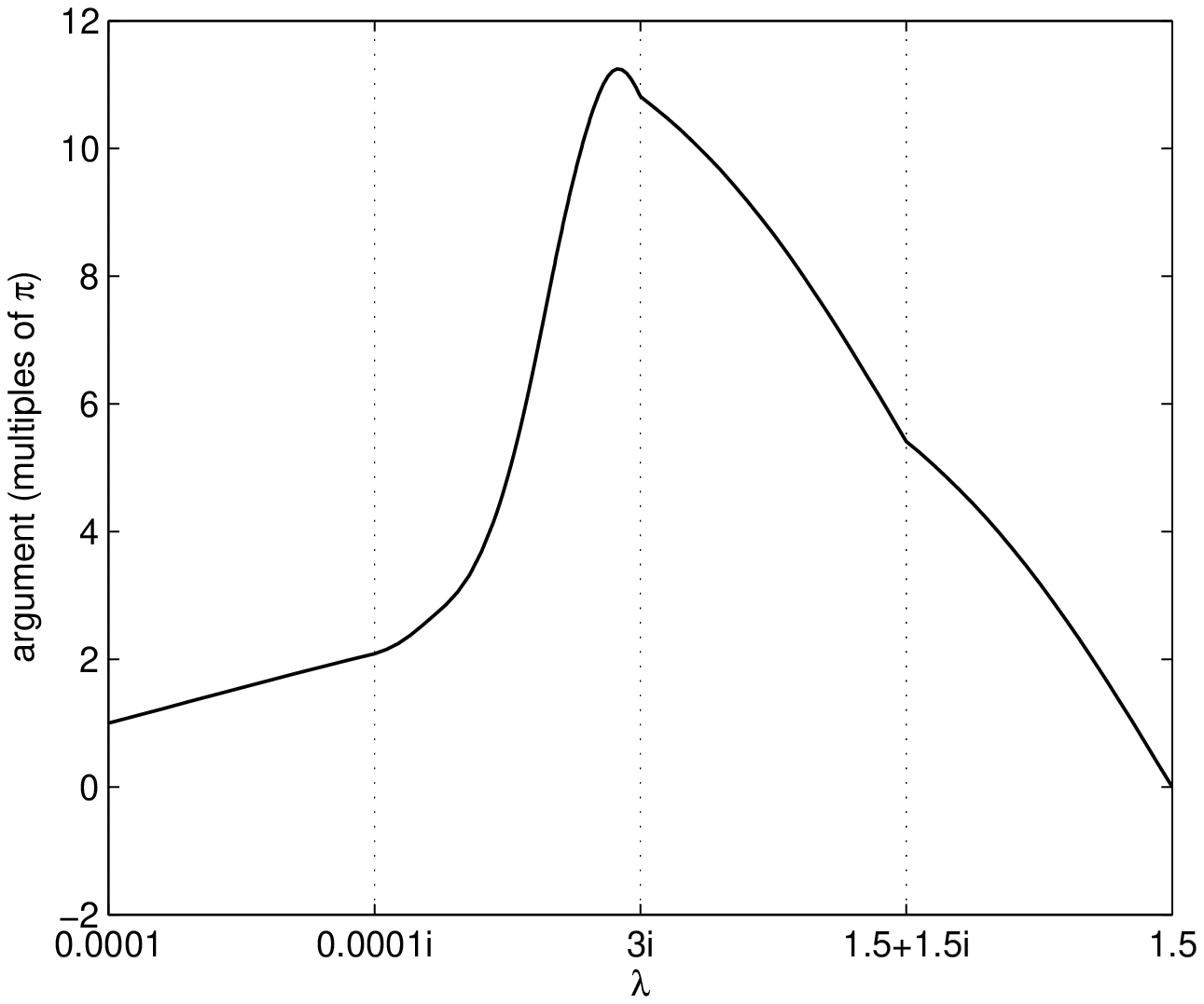}} 
\end{center}
\caption{The left figure shows the contour~$\cC$ in the complex plane.
  The right figure shows the argument of $D(\lambda)$ when $\lambda$
  transverses the top half of this contour (for the wrinkled front
  with $\delta=3$). The plot is split up in four parts, corresponding
  to the quarter circle, the segment along the imaginary axis, the
  diagonal segment, and the segment going down. The scale on the
  horizontal axis is linear on the first, third and fourth part and
  logarithmic on the second part.}
\label{contour}
\end{figure}

In our case, we have $c \approx 0.577$ and
$\|\mathrm{D}F(U_c)\|_{L^\infty} \approx 1.422$ for $\delta = 2.5$,
and $c \approx 0.548$ and $\|\mathrm{D}F(U_c)\|_{L^\infty} \approx
1.450$ for $\delta = 3$. So, the above estimate implies that any
eigenvalue~$\lambda$ satisfies $\mathrm{Re}(\lambda) < 1.5$ and
$\mathrm{Re}(\lambda)+|\mathrm{Im}(\lambda)| < 3$. Any unstable
eigenvalue is thus inside the contour bounded by these lines and the
imaginary axis; this contour is depicted in the left part of
Figure~\ref{contour} (the small semi-circle around the origin is to
exclude the double eigenvalue at the origin, as explained below). The
Evans function is analytic inside this contour, and thus we can count
the number of zeros inside the contour by computing the change of the
argument of the Evans function as $\lambda$ goes around the contour.

This suggests the following method. We divide the contour~$\cC$ in small
segments, compute the Evans functions at the end points of the
segments, find the change in the argument over each section and sum
these to find the total change over the whole contour. The number of
zeros is then the change of the argument divided by~$2\pi$. The
segments in which the contour is divided have to be chosen small enough
that the change of the argument over each segment is smaller
than~$\pi$. Adaptive procedures to achieve this have been proposed by,
amongst others, Ying and Katz~\cite{YK}, but for simplicity we
determined the partition of the contour by trial and error.

The Evans function, with our choice of initial conditions, satisfies
$D(\overline{\lambda}) = \overline{D(\lambda)}$. Hence, it suffices to
transverse half the contour~$\cC$. The plot at the right of
Figure~\ref{contour} shows how the argument of~$D(\lambda)$ changes as
$\lambda$ transverses the top half of the contour. We see that the
change over the whole contour is~$-2\pi$, and thus there is one zero
inside the contour. This is the zero around $\lambda = 0.0016$ shown
in Figure~\ref{evansw3qr}; there are no other unstable eigenvalues.

The operator~$\cL$ is invariant under translations and hence it has
two eigenvalues at $\lambda = 0$. Thus, the Evans function is zero at
$\lambda = 0$ and its argument is undefined. For this reason, the
contour~$\cC$ includes a small semi-circle of radius~$10^{-4}$ to
avoid the origin. It is of course possible that we miss some
eigenvalues because of this deformation. To guard against this, we
also compute the change in the argument of~$D(\lambda)$ as
$\lambda$~transverses a circle around the origin with
radius~$10^{-4}$. We find that the argument changes by~$4\pi$ and hence
the Evans function has two zeros inside the circle. This shows that we
have not missed any eigenvalues by the deformation of~$\cC$ around the
origin.

Another method is the parity-index argument 
(which we have not investigated here, but we include 
a description for completeness).
To determine if there are eigenvalues on the positive real axis, 
we could use that
\begin{equation*}
D(\lambda)\sim2\delta^{-K-1/2},
\end{equation*}
as $\lambda\to+\infty$, as we prove in Appendix~\ref{app:largelam}. 
A change in sign of the Evans function along the real axis indicates
an unstable eigenvalue. This can be exposed by comparing the 
sign of the Evans function in the asymptotic limit $\lambda\to+\infty$
to its sign elsewhere on the positive real axis, or the derivative
or second derivative of the Evans function along the real axis
evaluated at the origin. Note that to implement this 
technique we need to ensure the Evans function evaluated in
the neighbourhood of the origin and that generating the asymptotic
limit have been normalized, as far as their sign is concerned, 
in the same way. This may require keeping careful track of the 
the non-zero analytic multiples dropped
in our definition of the modified Evans function,
i.e.\ the normalization of the initial data,
the scalar exponential factor and the determinants
of $u_{-}$ and $u_{+}$
(which can be computed analytically in the asymptotic limit also).

\subsection{Convergence of the Evans function}\label{sec:convergence}
The results reported in Table~\ref{evanstable} show that the
values of~$\lambda$ for which the Evans function~$D(\lambda)$ is zero
converge as $K\to\infty$. However they do not show how the Evans
function scales for other values of $\lambda$ (not in the spectrum)
as $K$ becomes large. We see in Figure~\ref{evansp3} for the planar front,
that in the domain of interest with $\lambda$ of order $10^{-4}$, the Evans function
computed using the Drury--Oja approach with $K=24$, is of order $10^{-60}$. 
Similarly, for the wrinkled front, we see in Figures~\ref{evansw3qr} 
and~\ref{evansw2.5qr}, also computed using the Drury--Oja approach 
with $K=9$, the scale of the Evans function is of order $10^{-45}$ 
for $\lambda$ of order $10^{-4}$. These small results 
we find for the Evans function---far
smaller than machine precision---suggest that the computation might be
inaccurate because of round-off error.  This is not the case here. The
solution of the Drury--Oja equation~\eqref{eq:DO} does not contain
such small numbers. They only appear when the determinant is computed
at the matching point, which is a stable computation; see for example
Higham~\cite[\S13.5]{Higham}. The determinant is computed by
performing an LU-decomposition with partial pivoting and then
multiplying the diagonal entries. Both steps are numerically stable
even though the final result may be very small (as long as it does not
underflow). The factorization~\eqref{eq:factorization} for the planar front 
provides another view point: the Evans function for $K=24$ is the product of
$2K+1=49$ one-dimensional Evans functions, and if each of them is 
0.06---not a very small number---then their product is of the
order~$10^{-60}$, which is a very small number. 
 
For the modified Evans function we construct using the Riccati
approach, we see in Figure~\ref{evansw3ric}, that in the domain
of interest, it is of order $10^{-6}$ when $K=9$. If we reproduce
Figure~\ref{evansw3ric} for increasing values of $K$, we find 
that the modified Evans function actually grows.  
For large $K$, it increases roughly by a factor of $10$ per unit 
increase in $K$. However it is important to emphasize that
the solutions to the Riccati systems, that are used to construct
the modified Evans function, remain bounded. It is only
when we evaluate the determinant, that for example when $K=21$, 
the scale of the modified Evans function in  
Figure~\ref{evansw3ric} is of order $10^4$. 
Hence the question of the scaling of the Evans function comes
into play in the multi-dimensional context. 

In the one-dimensional case, the standard definition 
of the Evans function is~\eqref{defevans}
given in Alexander, Gardner and Jones~\cite{AGJ}. Of course
this definition stands modulo any analytic non-zero multiplicative
factor---afterall it's the zeros of the Evans function we are typically
after. A different appropriate normalization of the 
initial data $Y_0^\pm(\lambda)$ rescales the Evans function
by such a factor. Further it is common to drop the scalar
exponential factor. Hence for example the difference between
the Humpherys and Zumbrun Evans function~\eqref{eq:HZEvans}
and the standard Evans function is that the scalar
exponential factor is dropped and the initial data might
be normalized differently when it is $QR$-factorized.
For the modified Evans function in~\eqref{newdefevans}, 
we dropped the exponential factor and also the determinants
of some full rank submatrices in the determinant of the 
standard Evans function. 

However such analytic non-zero factors impact the 
scale of the Evans function in the multi-dimensional context,
as they depend on $K$. This can lead to growth or decay
of the Evans function employed away from the eigenvalues 
as $K$ is increased, depending on which factors are kept
and which are dropped. This is precisely what we have 
observed for the Drury--Oja and Riccati approaches,
which decay and grow, respectively, as $K$ increases.

Gesztesy, Latushkin and Zumbrun~\cite[\S4.1.4]{GLZ} show that the Evans
function, when suitably normalized, coincides with the 2-modified Fredholm
determinant of a certain operator related to $\cL$. Importantly, 
the 2-modified Fredholm determinant converges as $K\to\infty$.
They treat the self-adjoint case; the non-self-adjoint
operator $\cL$ we consider can, using suitable exponential weight functions, 
be transformed into a self-adjoint one (numerically though this
is not necessary and could introduce further complications due
to the introduction of exponential factors in the potential term).
Gesztesy, Latushkin and Zumbrun define the Evans function much like we have, 
through a Galerkin approximation, using the matching condition~\eqref{defevans}. 
Their normalization, natural to their context, involves ensuring
the Evans function does not depend on the coordinate system chosen and 
is invariant to similarity transformations (see 
Gesztesy, Latushkin and Makarov~\cite{GLM}).
Crucially, Gesztesy, Latushkin and Zumbrun demonstrate in 
Theorem~4.15 that their 2-modified Fredholm determinant
equals the Evans function multiplied by an exponential
factor (nonzero and analytic), whose exponent diverges as $K\to\infty$
(it can diverge to positive or negative infinity depending on the 
precise form of the coefficient matrix $A(x;\lambda)$).
Hence their Evans function diverges as well. Consequently, the 
2-modified Fredholm determinant is a natural candidate
for the canonical Evans function. 

With this in mind, our observation regarding the measured
divergence of the modified and Humpherys and Zumbrun Evans
functions in our experiments, is not too surprising. Their
divergence, the modified Evans function to large values,
and the Humpherys and Zumbrun Evans function to small values,
is simply a reflection of the over/under estimate of the 
crucial, correct multiplicative exponential factor that
equates these Evans functions with the 
2-modified Fredholm determinant of Gesztesy, Latushkin and Zumbrun,
which guarantees a finite limit as $K\to\infty$. 
It would be interesting to see whether one can explicitly compute 
the factor relating the Evans function of Gesztesy, Latushkin and
Makarov and the Evans function, as defined in this paper.

\section{Concluding remarks}\label{sec:conclu}
Our goal in this paper was to show, for the first time,
that spectral shooting methods can be extended to study
the stability of genuinely multi-dimensional travelling fronts.
It was already obvious to the people working in this field that this
can be done by a Fourier decomposition in the transverse dimension,
thus reducing the problem to a large one-dimensional problem; however,
no-one had yet managed to tackle the resulting large one-dimensional problems.  
This issue was overcome in theory by the work of Humpherys and
Zumbrun~\cite{HZ} and latterly that of Ledoux, Malham and
Th\"ummler~\cite{LMT}, which showed that enough information for
matching can be retained by integrating along the natural underlying 
Grassmannian manifolds (though these ideas have been 
well known in the control theory and quantum chemistry
for a long while---see for example Hermann and Martin~\cite{HM}
and Johnson~\cite{Jo}). The dimension of the Grassmannian
manifolds scale with the square of the order of 
the original problem. These methods brought the possibility
of multi-dimensional shooting within our grasp, but 
there were still challenges to be overcome. We have addressed
the most immediate of these in this paper:
\begin{enumerate}
\item How to project the problem onto the finite transverse
Fourier basis and ensure a numerically well-posed large system 
of linear ordinary differential spectral equations;
\item How to construct the genuinely multi-dimensional fronts, this
was a particular challenge which we managed to overcome
using the combined techniques of freezing
and parameter continuation (the latter for the unstable
multi-dimensional fronts);
\item Choose a non-trivial but well-known example problem, 
which exhibits the development of planar front instabilities 
(an important initial testing ground for us) but also two-dimensional
wrinkled fronts that themselves became unstable as a physical
parameter was varied;
\item Apply the Humpherys and Zumbrun continuous orthogonalization
as well as the Riccati methods to study the stability of the 
wrinkled fronts, in particular integrating a large system 
of order $196$---using $K=24$ transverse modes---hitherto
the largest system considered in Evans function calculations
was of order $6$ (see Allen and Bridges~\cite{AB}).
\item Ultimately demonstrating the 
\emph{feasibility}, \emph{accuracy}, \emph{numerical convergence} 
and \emph{robustness} of multi-dimensional shooting to determine eigenvalues.
\end{enumerate}
However there is still more to do. Though multi-dimensional shooting 
is competitive in terms of accuracy, how does it compare in terms
of computational cost? There are multiple sources of error including:
approximation of the underlying wrinkled front, domain truncation,
transverse Fourier subspace finite projection, approximation of
the solution to the projected one-dimensional spectral equations
and the matching position. Our experience is that the accuracy of the
underlying wrinkled front is especially important. 
Although we have endeavoured to keep all these errors small, not only
for additional confidence in our numerical results, but also for
optimal efficiency, we would like to know the relative influence of
each one of these sources of error in order to target our
computational effort more effectively.  Finally we would also like to
see how the method transfers to hyperbolic travelling waves and more
complicated mixed partial differential systems.

\section*{Acknowledgements}
SJAM attended a workshop at AIM in Palo Alto in May 2005 
on ``Stability criteria for multi-dimensional waves and patterns'' 
organised by Chris Jones, Yuri Latushkin, Bob Pego, Bj\"orn Sandstede
and Arnd Scheel that instigated this work.
The authors would like to thank Marcel Oliver for his 
invaluable advice and help at the beginning of this
project and also Gabriel Lord, Tom Bridges, Bj\"orn Sandstede
and Jacques Vanneste for useful discussions.
A large part of this work was conceived and 
computed while all four authors were visiting
the Isaac Newton Institute in the Spring of 2007.
We would like to thank Arieh Iserles and Ernst Hairer
for inviting us and providing so much support and 
enthusiasm. Lastly we would like to thank the 
anonymous referees whose helpful comments and
suggestions markedly improved the overall 
presentation of the paper.

\appendix

\section{Sectorial estimate on the spectrum}\label{app:sectorial}
We follow the standard arguments that prove that
the parabolic linear operator $\mathcal L$ is sectorial,
see for example Brin~\cite[p.~26]{Br}.
First consider the $L^2$ complex inner product with the spectral problem
\begin{equation*}
B\,\Delta U+c\partial_xU+\mathrm{D}F(U_c)\,U=\lambda U.
\end{equation*}
Premultiplying this spectral equation with $U^\dag$,
and integrating over the domain $\mathbb R\times\mathbb T$, we get
\begin{equation*}
-\int_{\mathbb R\times\mathbb T}\nabla U^\dag B\nabla U\,\mathrm{d}x\mathrm{d}y
+c\int_{\mathbb R\times\mathbb T}U^\dag\,\partial_xU
+U^\dag\,\mathrm{D}F(U_c)\,U\,\mathrm{d}x\mathrm{d}y
=\lambda\|U\|_{L^2}^2.
\end{equation*}
In this last equation, if we set $B=P^{\mathrm{T}}P$ where
$P\equiv\mathrm{diag}\bigl\{B_{11}^{1/2},~B_{22}^{1/2}\bigr\}$,
we get
\begin{equation}\label{eq:lambdaeqn}
\lambda\|U\|_{L^2}^2+\|P\nabla U\|_{L^2}^2=
c\int_{\mathbb R\times\mathbb T}U^\dag\,\partial_xU
+U^\dag\,\mathrm{D}F(U_c)\,U\,\mathrm{d}x\mathrm{d}y.
\end{equation}
If we consider the complex conjugate
transpose of the spectral equation, postmultiply with $U$, 
and integrate over the domain $\mathbb R\times\mathbb T$, we get
\begin{equation*}
\bar{\lambda}\|U\|_{L^2}^2+\|P\nabla U\|_{L^2}^2=
c\int_{\mathbb R\times\mathbb T}\partial_xU^\dag\,U
+U^\dag\,\bigl(\mathrm{D}F(U_c)\bigr)^{\mathrm{T}}\,U\,\mathrm{d}x\mathrm{d}y.
\end{equation*}
If we now add these last two equations then we get
\begin{align}
\mathrm{Re}(\lambda)\|U\|_{L^2}^2+\|P\nabla U\|_{L^2}^2=&\;
\int_{\mathbb R\times\mathbb T}U^\dag\,
\tfrac12\bigl(\mathrm{D}F+(\mathrm{D}F)^{\mathrm{T}}\bigr)(U_c)
\,U\,\mathrm{d}x\mathrm{d}y \nonumber\\
\leq&\;\|\mathrm{D}F(U_c)\|_{L^\infty}\|U\|_{L^2}^2.
\label{eq:lambdaeqnre}
\end{align}
Now dividing through by $\|U\|_{L^2}^2$ we get~\eqref{eq:sectorial1}

Next consider taking the imaginary part of \eqref{eq:lambdaeqn} followed by
the absolute value on both sides and using H\"older's inequality to get 
\begin{equation}\label{eq:lambdaeqnim}
|\mathrm{Im}(\lambda)|\,\|U\|_{L^2}^2\leq
c\|U\|_{L^2}\|\nabla U\|_{L^2}
+\|\mathrm{D}F(U_c)\|_{L^\infty}\|U\|_{L^2}^2.
\end{equation}
Adding inequalities \eqref{eq:lambdaeqnre} 
and \eqref{eq:lambdaeqnim} and using that
$\|P\nabla U\|_{L^2}^2\geq \kappa\,\|\nabla U\|_{L^2}^2$,
where $\kappa\equiv\min\{B_{11},B_{22}\}$, 
and the arithmetic-geometric mean inequality 
\begin{equation*}
c\|U\|_{L^2} \|\nabla U\|_{L^2}\leq \frac{c^2}{4\kappa}\|U\|_{L^2}^2
+\kappa\|\nabla U\|_{L^2}^2,
\end{equation*}
we get~\eqref{eq:sectorial2}---after dividing through by $\|U\|_{L^2}^2$.

\section{Evans function when spectral parameter is large}\label{app:largelam}
We follow the arguments for such results given in Alexander, Gardner
and Jones~\cite[p.~186]{AGJ} and Sandstede~\cite[p.~24]{Sand}.
We begin by rescaling $\tilde x=\sqrt{|\lambda|}\, x$ in \eqref{eq:bigsys}
for each $k=-K,-K+1,\ldots,K$ and taking the limit 
$|\lambda|\rightarrow\infty$ to get
\begin{equation*}
\partial_{\tilde x\tilde x} \hU_k 
- \mathrm{e}^{\mathrm{i}\,\mathrm{arg}\lambda}B^{-1} \hU_k =0.
\end{equation*}
Hence we have a decoupled system of $2K+1$ equations of the form
\begin{equation*}
\partial_{\tilde x} \begin{pmatrix} \hU \\ \hP \end{pmatrix}
= \begin{pmatrix} 
O_{2(2K+1)}               & I_2\otimes I_{2K+1} \\ 
\tilde B\otimes I_{2K+1}  & O_{2(2K+1)} 
\end{pmatrix}
\begin{pmatrix} \hU \\ \hP \end{pmatrix}
\end{equation*}
where $\tilde B=\mathrm{e}^{\mathrm{i}\,\mathrm{arg}\lambda}B^{-1}$.
This is a constant coefficient linear problem and the so-called
spatial eigenvalues $\mu$ are the roots of the characteristic polynomial
\begin{equation*}
\mathrm{det}\bigl((\mu^2 I-\tilde B)\otimes I_{2K+1}\bigr)
=\bigl(\mathrm{det}(\mu^2 I-\tilde B)\bigr)^{2K+1}.
\end{equation*}
where for any matrix $C$ we have used that 
$\mathrm{det}\bigl(C\otimes I_k\bigr)
\equiv\mathrm{det}\bigl(C^k\bigr)
\equiv\bigl(\mathrm{det}\,C\bigr)^k$.
Hence the spatial eigenvalues are $\pm\sqrt{B_{11}}$, 
$\pm\sqrt{B_{22}}$; each with algebraic multiplicity $2K+1$.
Associated with each of these four basic spatial eigenvalue
forms we have the $(2K+1)$-dimensional eigenspaces:
\begin{align*}
\mu=+\sqrt{B_{11}}\colon&\quad 
\mathrm{span}\{e_{2k}^{+}(+\sqrt{B_{11}})\colon k=-K,-K+1,\ldots,+K\},\\
\mu=-\sqrt{B_{11}}\colon&\quad 
\mathrm{span}\{e_{2k}^{-}(-\sqrt{B_{11}})\colon k=-K,-K+1,\ldots,+K\},\\
\mu=+\sqrt{B_{22}}\colon&\quad 
\mathrm{span}\{e_{2k+1}^{+}(+\sqrt{B_{22}})\colon k=-K,-K+1,\ldots,+K\},\\
\mu=-\sqrt{B_{22}}\colon&\quad 
\mathrm{span}\{e_{2k+1}^{-}(-\sqrt{B_{22}})\colon k=-K,-K+1,\ldots,+K\},
\end{align*}
where $e_{i}(\beta)$ is the vector with $1$ in position $i$, with
$\beta$ in position $2(2K+1)+i$, and zeros elsewhere. 
Hence as $|\lambda|\rightarrow\infty$, the Evans function 
has the asymptotic form 
\begin{align*}
D(\lambda)\sim&\;\mathrm{det}
\begin{pmatrix} 
I_2\otimes I_{2K+1} & -I_2\otimes I_{2K+1}\\
\tilde B^{1/2}\otimes I_{2K+1} &\tilde B^{1/2}\otimes I_{2K+1}
\end{pmatrix}\\
=&\;2\,\mathrm{det}\bigl(\tilde B^{1/2}\otimes I_{2K+1}\bigr)\\
=&\;2\,\bigl(\delta^{-1/2}
\mathrm{e}^{\mathrm{i}\,\mathrm{arg}\lambda}\bigr)^{2K+1}.
\end{align*}

\end{document}